\documentclass[final, fleqn]{elsarticle}

\usepackage[margin=2.54cm]{geometry}

\usepackage[english]{babel}
\usepackage[utf8x]{inputenc}
\usepackage{hyperref}
\usepackage{placeins}  
\usepackage{multirow}
\usepackage[normalem]{ulem}
\usepackage{amsmath,amssymb, empheq} 
\usepackage{lineno}
\numberwithin{equation}{section}  
\usepackage{amsfonts,nccmath,bbm,mathrsfs,mathtools}
\newtheorem{thm}{Theorem}
\newtheorem{lem}[thm]{Lemma}

\newdefinition{rmk}{Remark}
\newproof{pf}{Proof}
\newproof{pot}{Proof of Theorem \ref{thm2}}

\usepackage[table]{xcolor}
\usepackage{xcolor}
\colorlet{shadecolor}{orange!15}

\usepackage{float}  
\usepackage{algorithm}
\usepackage{algpseudocode}

\makeatletter
\newcommand*\dotp{\mathpalette\dotp@{.5}}
\newcommand*\dotp@[2]{\mathbin{\vcenter{\hbox{\scalebox{#2}{$\m@th#1\bullet$}}}}}
\makeatother

\newcommand{\bsym}[1]{\boldsymbol{#1}}
\newcommand{\norm}[1]{\lVert#1\rVert}

\newcommand{\sprod}[2]{\langle#1,#2\rangle}

\def\Diff{\text{Diff}(\mathcal{S})}

\definecolor{Green}{RGB}{5, 156, 16}

\definecolor{Bleu}{rgb}{0.3, 0.2, 1.0}

\definecolor{Both}{rgb}{0.55, 0.2, 0.65}

\nolinenumbers
\begin{document}

\begin{frontmatter}
\title{A projection-based Characteristic Mapping method for tracer transport on the sphere}
\author[myaddress]{\corref{corrauthor}Seth Taylor} 
\ead{seth.taylor@mail.mcgill.ca}
\author[myaddress]{Jean-Christophe Nave}
\address[myaddress]{\footnotesize \textit{Department of Mathematics and Statistics, McGill University, Montr\'{e}al, Qu\'{e}bec H3A 0B9, Canada}}
\cortext[corrauthor]{Corresponding author}

\begin{abstract}

A semi-Lagrangian Characteristic Mapping method for the solution of the tracer transport equations on the sphere is presented. The method solves for the solution operator of the equations by approximating the inverse of the diffeomorphism generated by a given velocity field. The evolution of any tracer and mass density can then be computed via pullback with this map. We present a spatial discretization of the manifold-valued map using a projection-based approach with spherical spline interpolation. The numerical scheme yields $C^1$ continuity for the map and global second-order accuracy for the solution of the tracer transport equations. Error estimates are provided and supported by convergence tests involving solid body rotation, moving vortices, deformational, and compressible flows. Additionally, we illustrate some features of computing the solution operator using a numerical mixing test and the transport of a fractal set in a complex flow environment.  

\end{abstract}

\end{frontmatter}
\nolinenumbers
\normalsize

\section{Introduction}

Numerical simulations of transport on the sphere are an indispensable tool in atmospheric and climate modeling. An efficient method of transporting multiple quantities, while preserving the scales present in the initial condition, is through composition with the inverse flow map generated by the advecting velocity field \cite{wiin1959application}. In this paper, we introduce a numerical framework, the Characteristic Mapping (CM) method, for the computation of this map in a spherical geometry.        

Existing numerical methods for the transport equation aim to compute the evolution of the advected quantity using the Eulerian, Lagrangian or semi-Lagrangian frameworks. An Eulerian scheme employs a spatial discretization on a static grid, allowing for ease of access to the solution and an easily parallelizable implementation. However, representing the solution on a static grid comes at the cost of time-stepping restrictions for explicit integration schemes imposed by the Courant-Friedrichs-Lewy condition. Implicit exponential integrators \cite{clancy2013use} and techniques for spatial adaptivity \cite{berger1989local, berger1984adaptive, flyer2010rotational} have been developed for Eulerian schemes. Lagrangian schemes alleviate these restrictions by following particle trajectories in the flow. The resulting solutions are generally less prone to numerical dissipation and have better stability properties \cite{staniforth1991semi}. A well-known difficulty arising in Lagrangian schemes is that an initially well-ordered set of points may become increasingly disordered over time \cite{welander1955studies,perlman1985accuracy}. Techniques to avoid excessive distortion include reseeding at regular time intervals \cite{barba2005advances, beale1985high, magni2012accurate, nordmark1991rezoning} and adaptive methods \cite{bergdorf2006lagrangian, bergdorf2005multilevel}. Another more recent approach, called indirect remeshing, interpolates the inverse flow map for the Lagrangian trajectories and then re-samples the initial advected quantity \cite{bosler2017lagrangian, bosler2013particle, bosler2014lagrangian}.  Finally, the semi-Lagrangian (SL) framework maintains a spatial discretization using an Eulerian grid while computing the evolution in the Lagrangian frame. Specifically, characteristics are traced back in time and the advected quantity is updated via interpolation. Related approaches include the Arbitrary Lagrangian-Eulerian method \cite{hirt1974arbitrary}, particle-mesh methods \cite{cotter2007remapped, cottet2014high}, flux-form \cite{lin1996multidimensional} and more recent hybrid schemes \cite{kaas2013hybrid, xiao2015rbf, shankar2018mesh}. Furthermore, conservative SL schemes have been developed and successfully applied for transport on the sphere \cite{lauritzen2010conservative, harris2011flux}. \par

Other popular SL schemes for the linear advection of sets and surfaces are Level-Set methods \cite{osher2004level}. Higher order accuracy may be generated within these methods by additionally transporting gradient information of the solution \cite{jetschemesNave, nave2010gradient,kohno2013new}. For quantities containing sub-grid features or for sets with poor regularity, the transport of gradient information may be infeasible. The CM method was developed to address this problem by Mercier et al.  \cite{mercier2020characteristic} for the linear transport of arbitrary sets in two and three dimensional periodic domains. The CM method provides a framework for the approximation of the inverse flow map generated by a velocity field. This map forms a path within the space of $C^1$ diffeomorphisms of the domain, allowing for a semi-discretization in time as a composition of submaps. Each of these submaps are computed using the semi-Lagrangian Gradient-Augmented Level Set (GALS) method \cite{nave2010gradient}, yielding an approximation of the entire backward trajectory map as a globally differentiable interpolant. These techniques have been further extended to solve the incompressible Euler equations by Yin et al. \cite{yin2021characteristic,yin2021characteristic3} in the same geometry. \par

The methods developed in \cite{mercier2020characteristic,yin2021characteristic,yin2021characteristic3} utilized an intrinsic spatial discretization of the submaps in an angular coordinate system. This in turn permitted the use of classical interpolation techniques relying upon an underlying vector space structure. The geometry of the sphere complicates this intrinsic approach due to the need of multiple overlapping coordinate charts. Here, we consider instead an extrinsic formulation for the spatial discretization of the backward characteristic map using $C^1$ spherical splines in a projection-based interpolation framework for manifold-valued data. The aim of this work is to present the design and analysis of this manifold-valued data approximation technique for the inverse flow map generated by a velocity field using the GALS and CM methods and demonstrate its utility for tracer transport on the sphere.  

The paper is organized as follows: section \ref{section:1} begins with the mathematical formulation of the CM method for the solution of the tracer transport equations. Section \ref{section:2} details the numerical implementation, outlining the spatial discretization via quadratic spline interpolation on the Powell-Sabin split along with the time evolution using the GALS method. Thereafter, we provide error estimates supported by convergence tests for four standard test cases involving a variety of flow environments in section \ref{results}. Finally, we conclude with a presentation of some unique features of the method by a zoom-in on the solution for the transport of a fractal set in a complex flow environment. \\

\section{Mathematical Formulation}
\label{section:1}
In this section we outline the mathematical formulation of the CM method for tracer transport on a surface $\mathcal{S}\subset \mathbb{R}^3$. We first describe the governing equations and solution strategy for the backward characteristic map as developed in \cite{mercier2020characteristic,yin2021characteristic}, adapted here for the geometry of $\mathcal{S}$. Thereafter, we demonstrate how this map may be associated with a solution operator for the transport equations.

\subsection{The Backward Characteristic Map}

Let $\Diff$ be the space of diffeomorphisms of $\mathcal{S}$. The central object of interest in the CM method is the family of diffeomorphisms $\bsym{X}_{[t,s]} \in \Diff$ generated by a velocity field $\bsym{u}: \mathcal{S}\times \mathbb{R} \to T\mathcal{S}$ over an interval of time $[t,s]$. We call these diffeomorphisms the characteristic maps as they provide a solution operator to the ordinary differential equation 

\begin{equation}
\label{integral_curves}
\dot{\bsym{\gamma}}(t) = \bsym{u}(\bsym{\gamma}(t),t) \, , \quad \bsym{\gamma}(s) = \bsym{x} \, ,
\end{equation}
for characteristic curves $\bsym{\gamma} : \mathbb{R} \to \mathcal{S}$ parameterized by the initial condition $\bsym{x} \in \mathcal{S}$. Fixing $s = t_0 < t$, we define the forward characteristic map $\bsym{X}_{[t_0,t]}: \mathcal{S} \times \mathbb{R} \to \mathcal{S}$ as the solution operator to \eqref{integral_curves} and the backward characteristic map $\bsym{X}_{[t,t_0]}$ as the inverse of $\bsym{X}_{[t_0,t]}$, that is 

\begin{equation}
\label{inverse_map}
\bsym{X}_{[t,t_0]} \circ \bsym{X}_{[t_0, t]} = \text{id}_{\mathcal{S}} \, .
\end{equation}The backward characteristic map provides a solution operator to \eqref{integral_curves} backwards in time. Differentiating the expression \eqref{inverse_map} with respect to time we get that $\bsym{X}_{[t,t_0]}$ satisfies an initial value problem of the form

\begin{equation}
\label{map_transport}
\begin{aligned}
\partial_t\bsym{X}_{[t,t_0]} + d\bsym{X}_{[t,t_0]}(\bsym{u}) &= 0 \, ,\\
\bsym{X}_{[t_0,t_0]} &= \text{id}_{\mathcal{S}} \, ,
\end{aligned}
\end{equation}where $d\bsym{X}_{[t,t_0]}: T\mathcal{S} \to T\mathcal{S}$ is the differential of the map. The CM method utilizes the semi-group structure of the characteristic map to compute the solution to \eqref{map_transport}. Consider partitioning an interval of time into $m$ subdivisions $[t_{i}, t_{i+1}] \subset [0, T]$ for $i \in \{0, 1, \dots, m-1\}$. We may decompose $\bsym{X}_{[T,0]}$ into a set of submaps $\bsym{X}_{[t_{i+1}, t_i]}: \mathcal{S}\to \mathcal{S}$ (defined in \cite{yin2021characteristic}), determined by the solutions of the following initial value problems
\begin{equation}
\label{eq:sub_map}
\partial_t \bsym{X}_{[t_{i+1}, t]}  + D\bsym{X}_{[t_{i+1}, t]}(\bsym{u})  = 0 \, , \quad \bsym{X}_{[t_{i+1}, t_{i+1}]} = \bsym{x} \, ,
\end{equation}backward in time to $t_{i}$. Using the semigroup property, the backward characteristic map to time $T$ is then obtained via composition of each submap, that is

\begin{equation}
\label{eq:global_map}
\bsym{X}_{[T,0]}(\bsym{x}) = \bsym{X}_{[t_1, 0]} \circ \bsym{X}_{[t_1, t_2]}  \dots \circ \bsym{X}_{[t_m, t_{m-1}]}(\bsym{x}).
\end{equation}We note that this formulation is geometrically intrinsic and holds without the need for a Riemannian metric structure on the surface. 

\subsection{The Solution Operator to the Transport Equations}

The equations governing the evolution of a tracer mixing ratio $\phi$ and fluid density $\rho$ transported by a given velocity field $\bsym{u}$ in the absence of sources or sinks are given by

\begin{subequations}
\label{transport_eqs}
\begin{align}
\partial_t \rho + \nabla_{\mathcal{S}} \cdot (\rho \bsym{u}) &= 0  \label{mass_transport}\, , \quad \rho(\bsym{x},0) = \rho_0 \, ,
\\
\partial_t (\rho \phi) + \nabla_{\mathcal{S}} \cdot (\rho \phi \bsym{u}) & = 0  \label{transport_eq}\, ,\quad \phi_0(\bsym{x},0) = \phi_0 \, , 
\end{align}
\end{subequations}
where $\nabla_{\mathcal{S}} \cdot $ is the surface divergence operator and $\phi_0, \rho_0$ are the initial conditions. Let $U_0 \subseteq \mathcal{S}$ be a fixed reference configuration in the fluid such that $U(t) = \bsym{X}_{[0,t]}(U_0)$ and let  $\bsym{X}_{[0,t]}(\bsym{x}) = \bsym{\alpha} \in U(t)$ be the Lagrangian particle position with position $\bsym{x} \in U_0$ in the Eulerian frame. In integral form, the continuity equations \eqref{transport_eqs} can be expressed as a coupled set of conservation laws 

\begin{subequations}
\label{conservation_law}
\begin{align}
 &\frac{d}{dt}\int_{U(t)} \rho(\bsym{\alpha},t) \mu(\bsym{\alpha}) =  0 \,,\label{mass_conservation}
 \\
 &\frac{d}{dt}\int_{U(t)} \rho(\bsym{\alpha},t) \phi(\bsym{\alpha},t) \mu(\bsym{\alpha}) =  0 \label{tracer_conservation}\, ,
 \end{align}
\end{subequations}
where $\mu \in \Omega^2(\mathcal{S})$ is a volume form. By a change of variables with the backward characteristic map, conservation of mass \eqref{mass_conservation} gives us that 

\begin{equation}
 \label{transport_theorem}
 \int_{U(t)} \rho(\bsym{\alpha},t)\mu(\bsym{\alpha}) = \int_{U_0} \rho_0(\bsym{x}) \mu(\bsym{x})  = \int_{U(t)} \bsym{X}_{[t,0]}^*(\rho_0 \mu) = \int_{U(t)} \rho_0 \circ \bsym{X}_{[t,0]} (\bsym{\alpha}) J_{\mu}(\bsym{X}_{[t,0]})(\bsym{\alpha}) \mu(\bsym{\alpha})\, .
 \end{equation}where $J_{\mu}(\bsym{X}_{[t,0]}) \in C^{\infty}(\mathcal{S})$ is the Jacobian determinant with respect to $\mu$ of the backward characteristic map. Since $\rho_0$ is well-defined and $U_0$ is arbitrary, we get that the solution to \eqref{mass_transport} is given by
 
\begin{equation}
\label{mass_solution_operator}
\rho(\bsym{x},t) = \left(\rho_0 \circ \bsym{X}_{[t,0]} \right) J_{\mu}(\bsym{X}_{[t,0]})(\bsym{x}) \, , \quad \forall \, (\bsym{x}, t) \in \mathcal{S} \times \mathbb{R}_+
\end{equation}Using a similar change of variables as in \eqref{transport_theorem} for the tracer density conservation \eqref{tracer_conservation} and applying \eqref{mass_solution_operator} we get that the evolution of the tracer mixing ratio is governed by

\begin{equation}
\label{solution_operator}
\phi(\bsym{x},t) = \phi_0 \circ \bsym{X}_{[t,0]}(\bsym{x})  \, , \quad \forall \, (\bsym{x}, t) \in \mathcal{S} \times \mathbb{R}_+ \, .
\end{equation}Using the relations \eqref{mass_solution_operator} and \eqref{solution_operator}, we see that $\bsym{X}_{[t,0]}$ defines a solution operator to the transport equations \eqref{transport_eqs} by pullback

\begin{equation}
\label{transport_solution_operator}
\bsym{X}^*_{[t,0]}:  \mathcal{F}(\mathcal{S}) \times \text{Dens}(\mathcal{S}) \to \mathcal{F}(\mathcal{S}) \times \text{Dens}(\mathcal{S}) \,, \quad 
(\phi_0, \rho_0) \mapsto (\phi_0 \circ \bsym{X}_{[t,0]}, \rho_0 \circ \bsym{X}_{[t,0]}  J_{\mu}(\bsym{X}_{[t,0]}) )
\end{equation}
where $\mathcal{F}(\mathcal{S})$ and $\text{Dens}(\mathcal{S})$ are the spaces of $\mathbb{R}$-valued functions and of densities on $\mathcal{S}$ respectively. 

\begin{rmk}
The pullback action defined by \eqref{transport_solution_operator} is the natural action of $\text{Diff}(\mathcal{S})$ on the product space $\mathcal{F}(\mathcal{S}) \times \text{Dens}(\mathcal{S})$. In general, pullback with the inverse flow map generated by the velocity field defines a solution operator for the Lie advection equation on $k$-forms.\\ \par
\end{rmk}

In the case that the velocity field is incompressible we have that $J_{\mu}(\bsym{X}_{[t,0]}) = 1$. Moreover, if the initial fluid density is also constant then the evolution of $\rho \phi$ is governed entirely by \eqref{solution_operator}. In the forthcoming numerical tests we will consider an initially constant density fluid and solve for the evolution of $\rho$ and $\phi$ using both compressible and incompressible velocity fields. \par

There are many benefits of solving for the entire solution operator to the transport equations \eqref{transport_eqs}, common to many SL schemes. A single numerical approximation of the map $\bsym{X}_{[t,0]}$ can be used to transport multiple quantities at the cost of the composition \eqref{solution_operator} \cite{wiin1959application}. Moreover, since the solution operator \eqref{transport_solution_operator} acts on the entire space $\mathcal{F}(\mathcal{S})$, the regularity of an advected quantity bears no constraint on its evolution \ref{last_section}. Since the transport of $(\phi, \rho)$ is defined through sampling with the numerically approximated flow map, the burden of spatial discretization is pushed away from the advected quantities and onto the inverse flow map. The ability to sample the map at arbitrary points in the domain, without having to recompute a trajectory for each sample point introduces an error due to both spatial and temporal discretizations \cite{wiin1959application}. The CM method provides a technique for the approximation of this map which strikes a balance between these two errors using an evolution via composition.\par

\section{Numerical Framework}
\label{section:2}

In this section we describe the numerical framework for the CM method for the tracer transport equations on the unit sphere centred at the origin, i.e. $\mathcal{S} = \mathbb{S}^2$, hereafter referred to as just the sphere. We note that the projection-based formulation of the method is not limited to the sphere and could be extended to a more general class of manifolds beyond the scope of what is presented here. We begin with a brief description of the solution algorithm in a generalized sense and remark upon the discretization of the backward characteristic map in the context of manifold-valued data approximation. \par 
Let $\mathcal{V}_h \subset C^1(\mathcal{S},\mathcal{S})$ be a finite dimensional interpolation space, defined with respect to a spatial discretization of $\mathcal{S}$ and denote $\mathcal{J}_h: C^1(\mathcal{S}, \mathcal{S}) \to \mathcal{V}_h$ as an interpolation operator projecting into this space. Let $t_n = n\Delta t$ and denote $\mathcal{X}_{[t_n,0]} \in \mathcal{V}_h$ as the numerical approximation of $\bsym{X}_{[t_n,0]}$. The CM method employs a backward semi-Lagrangian approach to compute $\mathcal{X}_{[t_{n+1},0]}$. Each iteration updates the map using a discrete analogue of \eqref{eq:sub_map} and \eqref{eq:global_map}:
 
 \begin{subequations}
\label{CMmethod}
\begin{align}
 \mathcal{X}_{[t_{n+1}, t_n]}(\bsym{x}) &= \mathcal{R}_{\Delta t}(\bsym{u}, \bsym{x})  \label{numerical_int}\, ,
 \\
\mathcal{X}_{[t_{n+1}, 0]}(\bsym{x}) &= \mathcal{J}_h [\mathcal{X}_{[t_n, 0]} \circ \mathcal{X}_{[t_{n+1}, t_n]}](\bsym{x}) \label{interpolation_step}\, ,
\end{align}
 \end{subequations}
which is initialized by setting $\mathcal{X}_{[0,0]}$ as the identity map. The value of the submap at $\bsym{x}$ is computed using a numerical integration scheme $\mathcal{R}_{\Delta t} : \mathfrak{X}(\mathcal{S} \times \mathbb{R}) \times \mathcal{S} \to \mathcal{S}$ providing an approximate solution to the ordinary differential equation 
\begin{equation}
\label{characteristics}
\dot{\bsym{\gamma}}(t)= \bsym{u}(\bsym{\gamma}(t),t)\, , \quad \bsym{\gamma}(t_{n+1}) = \bsym{x} \, ,
\end{equation}
backwards in time to $t_n$. We use the built-in interpolation framework of the GALS method \cite{nave2010gradient} to perform the interpolation step \eqref{interpolation_step}.  This method uses a local Hermite interpolation operator $\mathcal{J}_h$ to interpolate transported function values and derivative information which is supplied approximately by integrating for the footpoints of a compact $\varepsilon$-difference stencil \cite{nave2010gradient} about the vertices. \par

The use of classical interpolation techniques which rely upon an underlying vector space structure are complicated by the non-linear nature of the sphere and the space $C^1(\mathcal{S}, \mathcal{S})$. In the intrinsic approach to this approximation problem, the map is expressed in a local coordinate system where standard interpolation techniques and numerical integration schemes on linear spaces can be applied\footnote{We have implemented a method based on a tensor product of Hermite cubics and a higher-order Runge-Kutta integration scheme on a latitude-longitude spherical coordinate parametric space. The approach is simple and was observed to be highly accurate for many standard cases from \cite{nair2010class}, although limited in its scope due to the coordinate system.}. The need for multiple coordinate charts necessitates the use of overlapping meshes, such as the Yin-Yang grid \cite{kageyama2004yin}, along with the blending of interpolants over these grids. In contrast, an extrinsic approach to the problem can be taken by considering $\mathcal{S}$ as embedded in $\mathbb{R}^3$ and performing the interpolation or numerical integration scheme in the ambient Euclidean space \cite{giraldo1999trajectory}. The interpolant can be subsequently constrained to the surface through the use of a projection operator. In turn, the extrinsic formulation is void of artificial coordinate singularities and possesses the same order of accuracy as the underlying interpolation scheme in the ambient space.\par 

In this work, we consider a projection-based method for the spatial discretization of the backward characteristic map, similar in spirit to the projection-based interpolation methods for manifold-valued data \cite{gawlik2018embedding,grohs2019projection}. The sphere is discretized using a spherical triangulation allowing for the use of the higher-order accurate spherical spline interpolation methods outlined in \cite{lai2007spline,alfeld1996fitting}. The particular choice of spherical spline space considered consists of the quadratic spherical splines on the Powell-Sabin split. This space permits a globally $C^1$ interpolant from data consisting only of the function values and surface derivatives at the vertices. The resulting scheme offers computational efficiency, third-order accuracy, and algorithmic simplicity. The spatial discretization of the backward characteristic map is formed as a vector-valued spherical spline composed with a radial projection onto the sphere. Combining this interpolation technique with a higher-order numerical integration scheme results in a globally second-order accurate approximation of the backward characteristic map.

\subsection{Spatial Discretization}

We consider local spline spaces since they are constructed to be locally supported and globally differentiable; both desirable properties exploited by the GALS framework.  The particular choice of local spline space considered here are the quadratic spherical splines on the Powell-Sabin split. We note however that the use of other macro-element interpolation techniques such as those for cubic splines on the Clough-Tocher split or higher degree ($\geq 4$) polynomial splines are viable alternatives. A more comprehensive study of the GALS framework for more general spline spaces of higher degree and smoothness on other manifolds is the subject of current research. We begin by providing details on the particular method of interpolation along with error estimates describing the accuracy of the approximation.  

\subsubsection{Spherical Spline Interpolation on the Powell-Sabin Split}

Spherical analogues of bivariate spline spaces were introduced and studied in \cite{alfeld1996fitting, alfeld1996dimension,alfeld1996bernstein}. It was shown in \cite{alfeld1996fitting} that many of the bivariate macro-element techniques for spline interpolation on planar triangulations could be naturally carried over to the sphere using their construction. These techniques provide a powerful computational tool to perform local Hermite interpolation without the need to solve a linear system or construct an explicit basis. The sphere is discretized using a spherical triangulation $\mathcal{T}  = \{T_i\}_{i=1}^{N_{\Delta}}$ where $T_i \subset \mathcal{S}$, constructed by joining an associated set of vertices $\mathcal{V} = \{v_i\}_{i=1}^{N_v}$ along spherical arcs such that any two spherical triangles share at most one edge and their union covers the entire sphere. Let $\Pi_d$ be the space of homogeneous trivariate polynomials and denote $\Pi_d(\Omega)$ as the restriction of $\Pi_d$ to a subset $\Omega \subset \mathcal{S}$. The space of spherical splines of smoothness $r$ and degree $d$ on $\mathcal{T}$ is defined by \cite{alfeld1996bernstein}:

\begin{equation}
\label{ps_space}
S^r_d(\mathcal{T}) = \left\{s \in C^r(\mathcal{S}) \, : \,  s\rvert_{T} \in \Pi_d(T) \quad \forall \,T \in \mathcal{T} \right\}\, .
\end{equation} 
In order to obtain full approximation power for an interpolatory spline of degree $d < 3r + 2$, the technique of subdivision must be used \cite{lai2007spline}. In the case $d = 2$,  one can consider the Powell-Sabin (PS) split, originally introduced for bivariate spline interpolation on planar triangulations \cite{powell1977piecewise}. The subdivision is formed by joining the vertices and midpoints of each edge to the center of the triangle, resulting in six triangles associated with each macro-triangle $T$ (see left panel of figure \ref{fig:ps_split}). Let $\mathcal{T}_{PS}$ denote the PS refinement of $\mathcal{T}$ and let $\bsym{e}^1_i,\bsym{e}^2_{i} \in \mathbb{R}^3$ denote the normalized tangent vectors pointing along the edges connected to a vertex $\bsym{v}_i$. It was shown in \cite{alfeld1996fitting} that for every $f \in C^1(\mathcal{S})$ there exists a unique spline $s \in S^1_2(\mathcal{T}_{PS})$ satisfying the Hermite interpolation problem 

\begin{equation}
\label{Hermite_problem}
s(\bsym{v}_i) = f(\bsym{v}_i) \, ,\,\,\,  D_{\bsym{e}^1_i} s(\bsym{v}_i) = D_{\bsym{e}^1_i} f(\bsym{v}_i),  \,\,\,  D_{\bsym{e}^2_{i}} s(\bsym{v}_i) = D_{\bsym{e}^2_{i}} f(\bsym{v}_i)  \,, \quad  \forall \bsym{v}_i \in \mathcal{V} \, .
\end{equation}
 The linear functionals associated to the Hermite problem \eqref{Hermite_problem} on $\mathcal{T}_{PS}$ define an interpolation operator 
\begin{equation}
\label{Hermite_operator}
\mathcal{I}_h: C^1(\mathcal{S}) \to S^1_2(\mathcal{T}_{PS}).  
\end{equation}
The resulting spherical spline interpolant can be written explicitly using Bernstein-B{\'e}zier techniques, combining the $9$ pieces of data with the $C^1$ continuity conditions across the edges of the subdivision to map into the space $S_2^1(\mathcal{T}_{PS})$ as a set of $19$ coefficients for each macro-triangle $T$. Details on the explicit construction are given in the Appendix \ref{appendix}.
\par
Let $|T_i|$ be the maximum angular distance between any two points in $T_i$, and let $h = \max \{ |T_i|\}$. The interpolation operator \eqref{Hermite_operator} satisfies the following error estimate \cite{davydov2008interpolation,schumaker2015spline}:

\begin{equation}
\label{spline_error}
\norm{f - \mathcal{I}_h[f]}_{C(\mathcal{S})} \leq K h^{m+1} \norm{f}_{C^{m+1}(\mathcal{S})} \, , \quad \forall f \in C^{m+1}(\mathcal{S})
\end{equation}where $0 \leq m \leq 2$ and the constant $K$ depends on the smallest angle in the triangulation. Moreover, for any $\bsym{g} \in T\mathcal{S}$ the directional derivative of the interpolant satisfies

\begin{equation}
\label{spline_derivative_error}
\norm{D_{\bsym{g}}(f - \mathcal{I}_h[f])}_{C(\mathcal{S})} \leq K h^{m}\norm{f}_{C^{m+1}(\mathcal{S})} \, .
\end{equation}
\subsubsection{Projection-Based Spline Interpolation}
We now modify the spherical spline interpolation to interpolate a map $\bsym{X} : \mathcal{S} \to \mathcal{S}$. We consider the codomain of $\bsym{X}$ as embedded in $\mathbb{R}^3$ and interpolate the vector-valued map $\bsym{X} = (X_1, X_2, X_3)$, where $X_i: \mathcal{S} \to \mathbb{R}$ for each $i = 1,2,3$. The interpolation operator is defined as

\begin{equation}
\label{x_R}
\begin{aligned}
\mathcal{I}^{(3)}_{h} : \,  C^{1}({\mathcal{S} , \mathbb{R}^3}) &\to S^1_2(\mathcal{T}_{PS})^3
\\
\bsym{X} &\mapsto (\mathcal{I}_h(X_1), \mathcal{I}_h(X_2), \mathcal{I}_h(X_3)) \, .
\end{aligned}
\end{equation}
Since the interpolation operator is not constructed to enforce the constraint $\norm{\mathcal{I}_h^{(3)}[\bsym{X}](\bsym{x})} = 1$ globally, where $\norm{\,\cdot\,}$ is the Euclidean norm in $\mathbb{R}^3$, we compensate by composing the operator \eqref{x_R} with the radial projection onto the sphere

\begin{equation}
\label{proj}
\begin{aligned}
\mathcal{P} : C^{1}(\mathcal{S}, \mathbb{R}^3) \to  C^{1}(\mathcal{S}, \mathcal{S})\, , \,\,\, 
\bsym{X} \mapsto \bsym{X}/\norm{\bsym{X}} \, .
\end{aligned}
\end{equation}The interpolation operator used to define the spatial discretization of the backward characteristic map in \eqref{CMmethod} is thus given by

\begin{equation}
\label{Xtild}
\mathcal{J}_h : C^1(\mathcal{S}, \mathcal{S}) \to C^1(\mathcal{S}, \mathcal{S})\, , \quad \bsym{X} \mapsto \mathcal{P} \circ \mathcal{I}^{(3)}_h[\bsym{X}] \, .
\end{equation}Using this projection-based approach, the pointwise error is no worse than that of  $\mathcal{I}^{(3)}_h$, only up to a factor of 2 \cite{gawlik2018embedding}. This follows from an elementary estimate.  Since $\bsym{X}(\bsym{x}) \in \mathcal{S}$, we have that

\begin{equation}
\label{ptwise_error}
\begin{aligned}
\norm{\mathcal{J}_h(\bsym{X})(\bsym{x}) - \bsym{X}(\bsym{x})} &\leq \norm{\mathcal{P} \circ \mathcal{I}^{(3)}_h(\bsym{X})(\bsym{x}) - \mathcal{I}^{(3)}_h (\bsym{X})(\bsym{x})} + \norm{\bsym{X}(\bsym{x}) - \mathcal{I}^{(3)}_h(\bsym{X})(\bsym{x})}
\\
& \leq 2 \norm{\bsym{X}(\bsym{x}) - \mathcal{I}^{(3)}_h(\bsym{X})(\bsym{x})} \, .
\end{aligned}
\end{equation}
The differential of \eqref{Xtild} also satisfies a similar error estimate, but is no longer independent of the geometry of the sphere. We summarize a result on the approximation of the differential proven in \cite{hielscher2021approximating} for more general manifolds, here in the particular case of the sphere. Let $\bsym{\xi} = \mathcal{J}_h[\bsym{X}](\bsym{x})$ for some $\bsym{x} \in \mathcal{S}$ and suppose that $\bsym{v} \in T_{\bsym{x}}\mathcal{S}$. The differential of \eqref{Xtild} is given by

\begin{equation}
\label{d_Xtild}
d\mathcal{J}_h[\bsym{X}]_{\bsym{x}}: T_{\bsym{x}}\mathcal{S} \to T_{\bsym{\xi}}\mathcal{S} \,, \quad \bsym{v} \mapsto d\mathcal{P}_{\bsym{\xi}}\cdot D\mathcal{I}^{(3)}_h[\bsym{X}](\bsym{x}) \bsym{v} \, ,
\end{equation}
where $D\mathcal{I}_h^{(3)}[\bsym{X}](\bsym{x})$ is the Jacobian of \eqref{x_R} evaluated at $\bsym{x}$ using \eqref{SBB_derivative} and the differential of the projection at $\bsym{\xi} \in \mathbb{R}^3 \setminus \{0\}$ is given by

\begin{equation}
\label{proj_diff}
d\mathcal{P}_{\bsym{\xi}}: \mathbb{R}^3 \to T_{\bsym{\xi}}\mathcal{S} \, ,\, 
\bsym{v} \mapsto \frac{1}{\norm{\bsym{\xi}}}\left(\bsym{v} -  \frac{\sprod{\bsym{\xi}}{\bsym{v}}}{\norm{\bsym{\xi}}}\frac{\bsym{\xi}}{\norm{\bsym{\xi}}}\right) \,  .
\end{equation}The error in the differential of a map $F:\mathcal{S} \to \mathbb{R}^3$ is measured using the following norm

\begin{equation}
\norm{dF}_{2,\infty} \coloneqq \sup_{\bsym{x} \in \mathcal{S}} \norm{dF_{\bsym{x}}}_2
\end{equation}where $\norm{\,\cdot\,}_2$ is the matrix $2$-norm. As proven for more general manifolds in \cite{hielscher2021approximating}, if there is a positive constant $\epsilon < 1$ such that

\begin{equation}
\label{reach_condition}
\norm{\bsym{X}(\bsym{x}) - \mathcal{I}_h^{(3)}[\bsym{X}](\bsym{x})} \leq \epsilon \, , \quad \forall \bsym{x} \in \mathcal{S} \, ,
\end{equation}then the approximation of the differential of the map \eqref{Xtild} satisfies 
\begin{equation}
\label{differential_approx}
\norm{d\bsym{X} - d\mathcal{J}_h[\bsym{X}]}_{2, \infty}  \leq \norm{d\mathcal{I}^{(3)}_h[\bsym{X}] - d\bsym{X}}_{2, \infty} 
+ C(\epsilon)\left(\norm{d\mathcal{I}^{(3)}_h[\bsym{X}] - d\bsym{X}}_{2, \infty}  + \norm{d\bsym{X}}_{2, \infty}\right)\, ,
\end{equation}
where $C(\epsilon) = \epsilon(3 - 2\epsilon)(1-\epsilon)^{-1}$. 

\begin{rmk}
In the general setting of embedding-based interpolation methods for manifold-valued data, the pre-asymptotic behaviour of the approximation of the differential is dictated by the reach $\tau$ of the manifold, i.e. the size of the tubular neighbourhood about the embedded submanifold in which the projection is well-defined. A similar estimate to \eqref{differential_approx} holds with the condition that $\epsilon < \tau$ \cite{hielscher2021approximating}.
\end{rmk}

\subsubsection{Hermite Data Approximation}

The data supplied to the Hermite interpolation problem  \eqref{Hermite_problem} for the components of the map are given approximately using the technique of compact $\varepsilon$-difference stencilling about the vertices \cite{nave2010gradient}. Here, we give a description of this technique and an analysis of the effect of approximating the interpolation data \eqref{Hermite_problem} on the local truncation error using local tangent plane coordinates. We begin with a short description of this coordinate system and refer to \cite{davydov2008interpolation} for a more in-depth treatment on other embedded submanifolds.\par

Let $p \in \mathcal{S}$ and let $U_p \subset \mathcal{S}$ be an open set such that the projection onto the tangent plane at $p$ defined by $\pi_{p}: U_p \to T_p\mathcal{S}$ is invertible. Let $\bsym{p} \in \mathbb{R}^3$ also denote the normal vector to $\mathcal{S}$ at $p$ and let $\bsym{\gamma}_p^1, \bsym{\gamma}^2_p \in \mathbb{R}^3$ be unit vectors such that $\bsym{\gamma}_p^1 \times \bsym{\gamma}_p^2 = \bsym{p}$. Regarding $\mathcal{S}$ and $T_p\mathcal{S}$ as embedded in $\mathbb{R}^3$, then the tangent plane projection can be defined as

\begin{equation}
\label{tanplane_proj}
\pi_p(\bsym{q}) = \bsym{p} + \sprod{\bsym{q}}{\bsym{\gamma}^1_p}\bsym{\gamma}_p^1 + \sprod{\bsym{q}}{\bsym{\gamma}_p^2}\bsym{\gamma}^2_p  \,, \quad \forall q \in U_p\, .
\end{equation} Using a sufficient collection of charts $\varphi_p: U_p \to \mathbb{R}^2$ defined by

\begin{equation}
\varphi_p(q) = (\sprod{\bsym{q}}{\bsym{\gamma}^1_p}, \sprod{\bsym{q}}{\bsym{\gamma}_p^2}) : = (x_1, x_2) \, ,
\end{equation}
an atlas for $\mathcal{S}$ can be formed \cite{davydov2008interpolation}. We lastly note that the gradient of $\hat{f}_{p} = f \circ \varphi^{-1}_{p}$ at $\varphi_{p}(q) \in \varphi_{p}(U_{p})$, denoted by $\nabla\hat{f}_{p}(q)$, can be equivalently written as

\begin{equation}
\label{transition_grad}
 \nabla\hat{f}_{p}(q) = \nabla\hat{f}_{q}(q)J_{p, q}(q)  \, , 
\end{equation}where $J_{p, q}(q)$ is the Jacobian of the transition function $\varphi_{p} \circ \varphi_{q}^{-1}: \varphi_{q}(U_{q} \cap U_{p}) \to \varphi_{p}(U_{q} \cap U_{p})$ at $\varphi_{q}(q)$ given by  

\begin{equation}
\label{transition_jacob}
[J_{p, q}(q)]_{k,j}= \sprod{\bsym{\gamma}_{p}^{(k)}}{ \bsym{\gamma}_{q}^{(j)}} \, .
\end{equation} 
If we consider the chart $\varphi_{v_i}$ centred at the vertex $v_i \in \mathcal{S}$, the metric tensor reduces to the identity and the components of the surface gradient at $v_i$ are simply given by the partial derivatives of $\hat{f}$ evaluated at the origin. We approximate the function values and components of the surface gradient supplied in \eqref{Hermite_problem} by the $\mathcal{O}(\varepsilon^2)$ stencils

\begin{equation}
\label{eps_stencils}
\begin{aligned}
f(v_i) &\approx \frac{1}{4}\left(f_{(+,+)} + f_{(+,-)} + f_{(-,+)}+ f_{(-,-)}\right)
\\
\bsym{\gamma}^1_{v_i} \cdot \nabla_{\mathcal{S}}f(v_i) &\approx  \frac{1}{4 \varepsilon}  \left(f_{(+, -)} - f_{(-,-)} + f_{(+,+)} - f_{(-,+)} \right) \, ,
\\
\bsym{\gamma}^2_{v_i} \cdot \nabla_{\mathcal{S}}f(v_i) &\approx \frac{1}{4 \varepsilon} \left(f_{(-, +)} - f_{(-,-)} + f_{(+,+)} - f_{(+,-)} \right)  \, .
\end{aligned}
\end{equation}where $\varepsilon > 0$ is a constant and $f_{(\pm, \pm)} = f \circ \varphi_{v_i}^{-1}(\pm\varepsilon, \pm\varepsilon)$. We retain the local truncation error for \eqref{Hermite_problem} if we choose $\varepsilon$ such that the error incurred by \eqref{eps_stencils} is an $\mathcal{O}(h^3)$ approximation of function values and $\mathcal{O}(h^2)$ approximation of the derivatives. In practice the value of $\varepsilon$ is fixed beforehand such that the error incurred by \eqref{eps_stencils} is negligible compared to the local truncation error due to the interpolation for the discretizations considered. \par

A similar result regarding the accuracy of a PS interpolant for approximated Hermite data was given in \cite{davydov2008interpolation} using a different interpolation strategy and data fitting algorithm. We adapt their proof by transforming a spherical spline into a projected PS bivariate spline in local tangent plane coordinates. 
\begin{thm}
\label{approximate_Hermite_data}
Let $f \in C^3(\mathcal{S})$ and suppose that $D^{\bsym{\alpha}}f(v_i)$ is given approximately to order $\mathcal{O}(h^{3 - |\bsym{\alpha}|})$ for $|\bsym{\alpha}| \leq 1$. Let $s \in S^1_2(\mathcal{T}_{PS})$ be the solution to the interpolation problem \eqref{Hermite_problem} for the approximated data. Then the following holds

\begin{equation}
\label{approx_spline}
\begin{aligned}
\norm{f - s}_{C(\mathcal{S})} \leq K \left [ h^3 \norm{f}_{C^3(\mathcal{S})}  + \max_{v \in \mathcal{V}} \left( |f(v) - s(v)| + h \norm{\nabla f(v) - \nabla s(v)}\right)\right] \, ,
\end{aligned}
\end{equation}where $K$ depends on the smallest angle $\alpha$ in the triangulation and $\kappa = \sup_{p \in \mathcal{S}} \kappa(p)$ where

\begin{equation}
\kappa(p) = \max_{v \in \mathcal{V}_p} \left\{ \norm{J_{p, v}(v)}_2, \norm{J_{v, p}(v)}_{2}\right\} \, ,
\end{equation}for $\mathcal{V}_{p} \subset \mathcal{V}$ defined as the set of vertices contained in $U_{p}$. 
\end{thm}
\begin{pf}
Let $f_i, \bsym{\sigma}_{v_i}$ be the exact function values and directional derivatives at the vertices and let $\tilde{f}_i, \tilde{\bsym{\sigma}}_i$ be their approximations according to the statement of the theorem. Let $p \in S^1_2(\mathcal{T}_{PS})$ be the interpolatory spline solving \eqref{Hermite_problem} for the exact function and derivative data at the vertices and let $\bsym{v} \in \mathcal{S}$. Then using \eqref{spline_error} we have that 

\begin{equation}
\label{initial_estimate1}
|f(\bsym{v}) - s(\bsym{v})| \leq  Ch^3\norm{f}_{C^3(\mathcal{S})} + |p(\bsym{v}) - s(\bsym{v})| \, .
\end{equation}Suppose that $\bsym{v} \in \mathcal{S}$ is contained in the spherical triangle formed by the vertices $v_1, v_2, v_3$. The second term in \eqref{initial_estimate1} can be estimated by interpolating the functions $\hat{s} = s \circ \varphi_{\bsym{v}}^{-1}$ and $\hat{p} = p \circ \varphi_{\bsym{v}}^{-1}$ in a planar triangle $\hat{T}$ formed by the vertices $w_i = \varphi_{\bsym{v}}^{-1}(v_i)$ in the tangent plane at $\bsym{v}$. Let $\hat{h}$ the maximum edge length $\hat{T}$, which is bounded above by $h$. By the stability of the nodal minimal determining set for the space of bivariate splines on the PS split \cite{lai2007spline} and using \eqref{transition_grad} we get that 
\begin{equation}
\begin{aligned}
|s(\bsym{v}) - p(\bsym{v})| &= |\hat{s}(0,0) - \hat{p}(0,0)| \leq \max_{w_i} \left( |\hat{p}(w_i) - \hat{s}(w_i)| + \hat{h}\norm{\bsym{\sigma}_{v_i}J_{\bsym{v}, v_i}(w_i) - \tilde{\bsym{\sigma}}_{v_i}J_{\bsym{v}, v_i}(w_i)} \right)
\\
& \leq \max_{v_i} \left(|p(v_i) - s(v_i)| + h\kappa^2\norm{\bsym{\sigma}_{v_i} - \tilde{\bsym{\sigma}}_{v_i}} \right) \, ,
\\
\end{aligned}
\end{equation}which establishes the claim. \qed
\end{pf}
 
\subsection{Trajectory Computations}

As numerical integration scheme \eqref{numerical_int} for the footpoint computation we use the standard RK4 scheme in $\mathbb{R}^3$ backwards in time where the intermediate steps are projected back onto the sphere to evaluate the velocity field. Since the ordinary differential equation \eqref{integral_curves} is solved in Cartesian coordinates, without the use of a local coordinate system, numerical integration schemes which rely upon additive evolution without constraint will not remain on the sphere. Geometric numerical integration schemes have been devised to directly enforce this property  \cite{munthe1998runge,lewis2003geometric,iserles2000lie}. We have considered the fourth-order scheme RKMK4 \cite{munthe1998runge} and observed a similar accuracy to the approach taken here. In either case, these techniques do not rely upon an extension of the velocity field into the ambient space. Moreover, for $\Delta t$ small enough, the radial distance of the trajectory away from the sphere will be of the size of the local error of the integration scheme without projection \cite{hairer2006structure}. Consequently, the accuracy of the trajectory computation will not be compromised through the added projection computation. The extra computation of the normalization is however not required in the case that the velocity field is defined away from the sphere. 

\subsection{Algorithmic Implementation}

Given an implementation of the manifold-valued data approximation described by $\mathcal{J}_h$ and the numerical integration scheme described above, the computation of $\mathcal{X}_{[t,0]}$ can be summarized as follows. We precompute a spherical triangulation $\mathcal{T}$ and the initial stencils points by introducing a local orthonormal basis for the tangent plane at the vertices of the triangulation. The four stencil points $\varepsilon_{\pm,\pm}(v_i) \in \mathcal{S}$ in Cartesian coordinates at the vertex $v_i$ are given by $\varepsilon_{\pm,\pm}(v_i) = \pi_{v_i}^{-1}(\bsym{v}_i \pm \epsilon\bsym{\gamma}_{v_i}^1 \pm \epsilon\bsym{\gamma}_{v_i}^2, \bsym{v}_i \pm \epsilon\bsym{\gamma}_{v_i}^2  \pm \epsilon\bsym{\gamma}_{v_i}^1)$ where $\pi^{-1}_{v_i}$ is the inverse of the tangent plane projection \eqref{tanplane_proj} at $v_i$. We then apply the approach described by \eqref{CMmethod}: at each time step, the footpoints of the four stencil points about the vertices are computed by numerical integration backwards in time. The map from the previous iteration, given as an interpolant, is then evaluated at these foot points. The resulting values are then supplied to the compact stencilling \eqref{eps_stencils}, providing the data for the interpolation \eqref{Hermite_problem} in each component of the map. Finally, a new interpolant is formed as in \eqref{interpolation_step} by projecting back onto the space $\mathcal{P}(S^1_2(\mathcal{T}_{PS})^3)$ using $\mathcal{J}_h$, defining the next submap. \par

\subsubsection{Submap Decomposition}

The compositional structure of the evolution of the inverse flow map permits the use of the technique of submap decomposition to improve the accuracy of the method. This technique has been successfully applied to resolve the fast-growing gradients of the vorticity in the incompressible Euler equations \cite{yin2021characteristic}. The remapping strategy utilizes the fact that the initial condition of \eqref{map_transport} is given by the identity map. Since the identity map need not be interpolated, the footpoint computations on the first iteration incur no error due to spatial discretization. This property can be utilized along with the semi-group structure of the map to garner better accuracy in the solution. At time steps $\tau_i \in \{t_0, t_1, \dots, t_n\}$ where $0 < i \leq n-1$ we store the map $\mathcal{X}_{[\tau_i, 0]}$ in memory and reinitialize $\mathcal{X}_{[\tau_i, \tau_i]} = \text{id}_{\mathcal{S}}$. A submap is then computed over the interval $[\tau_{i}, \tau_{i+1}]$ using \eqref{CMmethod} on the time steps $t_j$ such that $\tau_i \leq t_j < \tau_{i+1}$. If the remapping is performed $m$ times then the numerical approximation of the backward characteristic map at $t_n$ is computed using the semi-group property as

\begin{equation}
\label{remapping}
\mathcal{X}_{[t_n,0]} = \mathcal{X}_{[\tau_1, 0]} \circ \mathcal{X}_{[\tau_2, \tau_1]} \circ \dots \circ \mathcal{X}_{[t_n, \tau_m]} \,.
\end{equation} 
Since the submaps forming the decomposition \eqref{remapping} are all initialized as the identity map, at each remapping step we effectively reset the local truncation error due to spatial discretization accumulated over the previous iterations. The computation of the backward characteristic map using \eqref{remapping} thus improves the accuracy of the method at the expense of increased memory allocation. \par

\begin{rmk}
The result of the approximation \eqref{remapping} is no longer an element of $\mathcal{P}(S_2^1(\mathcal{T}_{PS})^3)$ but instead in the $m$-times composed polynomial space $\mathcal{P}(S_2^1(\mathcal{T}_{PS})^3) \circ \dots \circ \mathcal{P}(S_2^1(\mathcal{T}_{PS})^3)$. The technique of submap decomposition approximates the entire map to time $t_n$ in a space with a polynomial degree of approximation which increases with the number of compositions of fixed polynomial degree forming \eqref{remapping}. 
\end{rmk}

\begin{rmk}
The remapping steps can be statically enforced or adaptively determined throughout the computation. The adaptive strategy employed in \cite{yin2021characteristic} for incompressible flows utilized the error in the Jacobian determinant of the map for some prescribed tolerance. A more ubiquitous adaptive strategy, independent of the compressibility of the velocity field, could be designed based on the conservation laws \eqref{conservation_law}. Although, this would require an additional numerical quadrature scheme, introducing a possibly costly computation at each iteration. A low-cost alternative adaptive remapping strategy is to integrate a set of Lagrangian particles forwards in time using a high-order integration scheme. A remapping step can be then be initialized using a prescribed tolerance on the error between the backward characteristic map evaluated at the location of these passive particles and the identity map. We refer to \cite{mercier2020characteristic} for further details on this strategy. In general, the most effective adaptive remapping strategies will be problem specific, depending on the particular properties of the advected quantity and transporting velocity field.
\end{rmk}

\subsection{Conservation Properties}
\label{sec:conservation_properties}
As a consequence of computing the evolution of the tracer and mass densities through pullback, the method possesses global mass conservation. By the change of variables formula with the $C^1(\mathcal{S}, \mathcal{S})$ approximation of the map we have that  

\begin{equation}
\label{mass_integral}
\int_{\mathcal{S}} \rho_0 \mu = \int_{\mathcal{S}} \mathcal{X}_{[t,0]}^*(\rho_0\mu) \,,
\end{equation}which holds similarly for the tracer density $\rho \phi$. A similar conservation property holds in a more local sense over sets $E \subset \mathcal{S}$, but involves an integration of $\mathcal{X}_{[t,0]}^*(\rho_0\mu)$ over the region $\mathcal{X}_{[0,t]}(E)$ where $\mathcal{X}_{[0,t]}$ is the spatial inverse of the numerically approximated inverse flow map. In the numerical results we have opted for an evaluation of the global mass conservation and the pointwise error in the density. These two quantities are computed as follows. Let the volume form $\mu$ be the two-form resulting from the interior product of the normal vector field with the Euclidean volume restricted to the sphere, that is

\begin{equation}
\mu(\bsym{x}) = x dy\wedge dz + y dz \wedge dx + z dx \wedge dy = \iota_{\bsym{n}}(dx \wedge dy \wedge dz)(\bsym{x}) \quad \forall \bsym{x} \in \mathcal{S} \,,
\end{equation}where $\iota_{\bsym{n}}$ is the interior product with the normal vector field $\bsym{n}(p) = \bsym{p} \in \mathbb{R}^3$. Consider the spherical coordinate chart defined by
\begin{equation}
\label{spherical_coords}
\eta^{-1} : [0,2\pi) \times (0,\pi) \to \mathcal{S} \subset \mathbb{R}^3 \,, \quad (\lambda, \theta) \mapsto (\sin\theta\cos\lambda, \sin\theta\sin\lambda, \cos\theta) \, .
\end{equation}In spherical coordinates, we get that the mass integral becomes 

\begin{equation}
\label{mass_integral_sph_coords}
\begin{aligned}
\int_{\mathcal{S}} \mathcal{X}_{[t,0]}^*\mu  &= \int_0^{\pi} \int_{0}^{2\pi}(\mathcal{X}_{[t,0]}\circ \eta^{-1})^*\mu 
\\
&= \int_0^{\pi} \int_{0}^{2\pi}\bigg[\sum_{\bsym{i} \in \sigma(1,2,3)}\mathcal{X}_{[t,0]}^{i_1}\circ \eta^{-1} \cdot  d(\mathcal{X}^{i_2}_{[t,0]}\circ \eta^{-1}) \wedge d(\mathcal{X}_{[t,0]}^{i_3}\circ \eta^{-1})\bigg](\lambda, \theta) d\lambda d\theta \,,
\end{aligned}
\end{equation}using the property $(\eta^{-1})^*\mathcal{X}_{[t,0]}^* = (\mathcal{X}_{[t,0]} \circ \eta^{-1})^*$ and where $\sigma(1,2,3)$ is the set of cyclic permutations of $\{1,2,3\}$ and the differentials $d\mathcal{X}^i_{[t,0]}\circ \eta^{-1}_T$ are computed via chain rule. We use a local coordinate expression for the Jacobian determinant for the evaluation of the pointwise error. Let $\bsym{\xi} = \mathcal{X}_{[t,0]}(\bsym{x}) \in \mathcal{S}$ an introduce the orthonormal bases $\{\gamma^{(1)}_{\bsym{x}}, \gamma^{(2)}_{\bsym{x}}\}$ and $\{\gamma^{(1)}_{\bsym{\xi}}, \gamma^{(2)}_{\bsym{\xi}}\}$ where for $T_{\bsym{x}}\mathcal{S}$ and $T_{\bsym{\xi}}\mathcal{S}$ respectively. The Jacobian determinant can be written as

\begin{equation}
J_{\mu}(\mathcal{X}_{[t,0]})(\bsym{x}) = \text{det}(d(\varphi_{\xi} \circ \mathcal{X}_{[t,0]} \circ \varphi_x^{-1}))(x_1, x_2) = \text{det}(d\hat{\mathcal{X}}_{[t,0]})(x_1, x_2)\,.
\end{equation}Since the charts are centred at $\bsym{x}$ and $\bsym{\xi}$ respectively, the components of the Jacobian matrix are simply given by

\begin{equation}
\left[d\hat{\mathcal{X}}_{[t,0]}\right]_{i,j} = \sprod{\gamma_{\bsym{\xi}}^j}{\nabla \mathcal{X}_{[t,0]}\cdot \gamma_{\bsym{x}}^{i}} 
\end{equation}where $\nabla \mathcal{X}_{[t,0]}$ is the Euclidean gradient of the map computed using \eqref{d_Xtild}.

\subsection{Error Estimates}
\label{sec:error_estimates}
In this section we provide error estimates for the solution of \eqref{map_transport} and \eqref{transport_eqs} using the CM method, yielding a theoretical justification of the results in section \ref{results}. We will use the following quantity to measure the error between two continuous maps $F, G:\mathcal{S} \to \mathcal{S}$ 

\begin{equation}
\label{CSS_norm}
\norm{F - G}_{C(\mathcal{S}, \mathcal{S})} := \sup_{x \in \mathcal{S}} d_{\mathcal{S}}(F(\bsym{x}), G(\bsym{x})) \, ,
\end{equation}which is the supremum of the arc length distance on the sphere between $F(\bsym{x})$ and $G(\bsym{x})$. Combining the estimates \eqref{ptwise_error}, \eqref{differential_approx} with \eqref{approx_spline} for each component of $\bsym{X}_{[t,0]}$ yields an estimate on the local truncation error due to the spatial discretization of the backward characteristic map. In particular, if the data supplied to \eqref{Hermite_problem} satisfy the conditions of theorem \eqref{approximate_Hermite_data} and \eqref{reach_condition} holds, then for a fixed $t$ there exists a $C > 0$ such that 

\begin{equation}
\label{spatial_error_map}
\norm{\mathcal{J}_h[\bsym{X}_{[t,0]}] -\bsym{X}_{[t,0]}}_{C(\mathcal{S},\mathcal{S})} \leq  Ch^3 \, .
\end{equation}\begin{lem}
Suppose that $T \in C^1(\mathcal{S}, \mathcal{S})$ satisfies $\norm{T - \text{id}_{\mathcal{S}}}_{C(\mathcal{S}, \mathcal{S})} = \mathcal{O}(\Delta t)$, then for all $f : \mathcal{S} \to \mathbb{R}$ sufficiently differentiable the following holds:

\begin{equation}
\label{bruce_lemma}
\norm{\mathcal{I}_h[f \circ T] - \mathcal{I}_h[f] \circ T}_{C(\mathcal{S}, \mathcal{S})} = \mathcal{O}(\Delta t h^2) \, .
\end{equation}\end{lem}
\begin{pf}
Using \eqref{spline_error} an analogous estimate was given in \cite{yin2020diffusion} (Lemma 3.1). 
\end{pf}

\begin{thm}
Using an $s$-stage RK integration scheme and quadratic spherical spline interpolant, the approximation $\mathcal{X}_{[T,0]}(\bsym{x})$ is consistent with $\bsym{X}_{[T,0]}(\bsym{x})$ to order

\begin{equation}
\label{CM_method_error}
\norm{\bsym{X}_{[T,0]} - \mathcal{X}_{[T,0]}}_{C(\mathcal{S},\mathcal{S})} = \mathcal{O}(T\min(h^2, h^3\Delta t^{-1}) + T\Delta t^s) \, .
\end{equation}
\end{thm}

\begin{pf}
Suppose initially that $\mathcal{X}_{[t_{n-1},0]} = \bsym{X}_{[t_{n-1},0]}$. The inverse of the submap $\mathcal{X}_{[t_n, t_{n-1}]}$ computed using an $s$-stage RK integration scheme forwards in time is an $\mathcal{O}(\Delta t^{s+1})$ approximation to the inverse of $\bsym{X}_{[t_n,t_{n-1}]}$ as measured using \eqref{CSS_norm}. Letting $\mathcal{E}_{[t_n, t_{n-1}]} = \bsym{X}_{[t_n,t_{n-1}]} \circ \mathcal{X}^{-1}_{[t_n, t_{n-1}]}$ we split the error incurred by \eqref{CMmethod} as follows

\begin{equation}
\label{initial_estimate}
\begin{aligned}
\norm{\mathcal{X}_{[t_n,0]}- \bsym{X}_{[t_n,0]}}_{C(\mathcal{S}, \mathcal{S})}  &\leq \norm{\bsym{X}_{[t_{n-1},0]} \circ \mathcal{X}_{[t_n, t_{n-1}]}- \bsym{X}_{[t_n,0]}(\bsym{x})}_{C(\mathcal{S}, \mathcal{S})}\\
&+\norm{\mathcal{J}_h[\bsym{X}_{[t_{n-1},0]} \circ \mathcal{X}_{[t_n, t_{n-1}]}] - \bsym{X}_{[t_{n-1},0]} \circ \mathcal{X}_{[t_n, t_{n-1}]}}_{C(\mathcal{S}, \mathcal{S})}
\\
& \coloneqq E_{\Delta t} +  E_{h,\Delta t}\, .
\end{aligned}
\end{equation}Using the Lipschitz continuity of $\bsym{X}_{[t_{n-1},0]}$ we can estimate the first term as 

\begin{equation}
\begin{aligned}
E_{\Delta t} &= \norm{(\bsym{X}_{[t_{n-1}, 0]} - \bsym{X}_{[t_{n-1}, 0]} \circ \mathcal{E}_{[t_{n-1}, t_n]} ) \circ \mathcal{X}_{[t_n, t_{n-1}]}}_{C(\mathcal{S}, \mathcal{S})}
\\
&\leq \norm{\bsym{X}_{[t_{n-1}, 0]} - \bsym{X}_{[t_{n-1}, 0]} \circ \mathcal{E}_{[t_{n}, t_{n-1}]}}_{C(\mathcal{S}, \mathcal{S})} 
\\
&\leq C \norm{\text{id}_{\mathcal{S}} -\mathcal{E}_{[t_n, t_{n-1}]}}_{C(\mathcal{S}, \mathcal{S})}= \mathcal{O}(\Delta t^{s+1}) \, ,
\end{aligned}
\end{equation}where $C$ is the Lipschitz constant of $\bsym{X}_{[t_{n-1},0]}$. The second term in \eqref{initial_estimate} can be estimated by applying \eqref{spatial_error_map} and \eqref{bruce_lemma}, yielding

\begin{equation}
\begin{aligned}
 E_{h, \Delta t} &\leq \norm{\mathcal{J}_h [\bsym{X}_{[t_{n-1},0]} \circ \mathcal{X}_{[t_n, t_{n-1}]}] - \mathcal{J}_h [\bsym{X}_{[t_{n-1},0]}] \circ \mathcal{X}_{[t_n, t_{n-1}]}}_{C(\mathcal{S}, \mathcal{S})} 
 \\
 &+ \norm{(\mathcal{J}_h [\bsym{X}_{[t_{n-1},0]}]  - \bsym{X}_{[t_{n-1},0]} )\circ \mathcal{X}_{[t_n, t_{n-1}]}}_{C(\mathcal{S}, \mathcal{S})} 
 \\
 &= \mathcal{O}(\Delta t h^2 + h^3) \, .
\end{aligned}
\end{equation}
The global truncation error for $N_t = T \Delta t^{-1}$ steps is then given by \eqref{CM_method_error}, establishing the claim.\qed
\end{pf}
The order of convergence of the approximation of \eqref{solution_operator} using $\mathcal{X}_{[t,0]}$ will in general by dictated by the modulus of continuity of $\phi_0$ when measured in the sup-norm. It follows that for any Lipschitz continuous initial condition $\phi_0$ the order of convergence of the solution will be equivalent to the order given by \eqref{CM_method_error}. In the weaker case that $\phi_0$ is only of bounded variation, then the convergence of the approximated solution will be of the same order as \eqref{CM_method_error} in the $L^1(\mathcal{S})$ norm. This follows from the fact that functions of bounded variation are almost everywhere differentiable and hence Lipschitz continuous up to a set of measure zero. \par

We can give also a uniform bound on the pointwise error in the density in terms of the error in the Jacobian determinant of the map. If we assume that the error in the differential of the map is small, then the sup-norm error of the approximated density satisfies

\begin{equation}
\begin{aligned}
\norm{J_{\mu}(\mathcal{X}_{[t,0]})-J_{\mu}(\bsym{X}_{[t,0]}) }_{L^{\infty}(\mathcal{S})} & \approx \norm{\text{tr} \left((d\bsym{X}_{[t,0]})^{-1}(d\mathcal{X}_{[t,0]} - d\bsym{X}_{[t,0]})\right)}_{L^{\infty}(\mathcal{S})}  
\\
&= \mathcal{O}(\norm{d\mathcal{X}_{[t,0]} - d\bsym{X}_{[t,0]}}_{2,\infty}) \, .
\end{aligned}
\end{equation}

\section{Numerical Tests} \label{results}

In this section we present numerical tests of the CM method for the solution of \eqref{transport_eqs} and \eqref{map_transport}. We begin with a numerical verification of the method and our implementation using a standard test case suite as outlined in \cite{nair2010class, nair2008moving}. The tests performed affirm the theoretically expected rate of convergence for the characteristic map \eqref{CM_method_error} along with the tracer and density errors in each case. Thereafter, we consider tests which illustrate certain properties of the advective nature of the error introduced by the method. These tests include the preservation of correlations between multiple advected quantities, the conservation of mass/area, and a simulation of the evolution of a fractal set in a complex flow environment.

\subsection{Numerical Verification}
\label{convergence_tests}

In this section we present a numerical verification of the CM method and our implementation using convergence tests based on a standard suite of velocity fields \cite{nair2008moving,nair2010class}. The tests are organized into solid body rotation \ref{sec:test_sbr}, deformational \ref{sec:deform_flow}, and compressible \ref{sec:test_compressible} velocity fields. In supplement to these tests we demonstrate how the CM method exactly preserves non-linear correlations between multiple advected quantities.

\subsubsection{Implementation Details}
\label{additional_details}
The numerical tests were all implemented in Python and run on a Linux workstation with an Intel core i5-8250U (8 logical processors) with 16 GB of RAM. The spherical triangulation is constructed using the Python package Stripy  \cite{moresi2019stripy} which provides a wrapper to the package STRIPACK \cite{renka1997algorithm}. The point in triangle querying was performed using the Python binding supported by the package Libigl \cite{libigl}. The code was written in a high-level language and the tests were performed on a laptop computer -- it is not an optimized implementation of the method. The run times are however modest for our purposes, taking approximately $8$ minutes for $100$ time steps with $N_{\Delta} = 327680$ for instance. It is observed that the method scales linearly with the number of triangles and the number of time steps $N_t$, i.e. the computational time is  $\mathcal{O}(N_{\Delta} + N_t)$, since explicit time stepping is used. The dominant contribution to the computational time is the evaluation of $\mathcal{X}_{[t_n,0]}$ at the footpoints of the $\epsilon$-difference stencils, computed at each iteration of the method, requiring the point in triangle querying. Quadtree data structures on the triangulation, along with parallelization of the footpoints calculations, could be implemented to improve the overall performance of the method.

The tests are all performed on successive uniform refinements of an icosahedral discretization of the sphere \cite{baumgardner1985icosahedral} (see table \ref{ico_table}). We note however that the formulation and implementation are essentially agnostic to the particular spherical triangulation. We demonstrate the convergence rate \eqref{CM_method_error} by refining $\Delta t = T/N_t$ using $N_t = 2^k + 10$ where $k$ is the number of refinements of the icosahedral discretization. The proportionality between $h$ and $\Delta t$ is chosen to reflect the global error coming from both the temporal and spatial discretizations. The chosen refinement is not a requirement; as a backward semi-Lagrangian evolution technique, the method does not possess time stepping restrictions based on the grid size. A value of $\varepsilon = 10^{-5}$ is chosen for the $\varepsilon$-difference stencils \eqref{eps_stencils} in each test which effectively limits the machine precision to approximately $10^{-12}$. \par

\begin{table}
\centering
\begin{tabular}{|c||c|c|c|c|c|c|c|c|c|}
\hline
$k$ & 0 & 1 & 2 & 3 & 4 & 5 & 6 & 7 & 8 \\
\hline\hline 
$N_v$  & 12 & 42 & 162 & 642 & 2562 & 10242 & 40062 & 163842 & 655362\\
\hline
$N_{\Delta}$ & 20 & 80 & 320 & 1280 & 5120 & 20480 & 81920 & 327680 & 1310720 \\
\hline
$h$ & 1.10715 & 0.62832 & 0.32637 & 0.16483 & 0.08263 & 0.04134 & 0.02067 & 0.01034 & 0.00517 \\ \hline 
\end{tabular}
\caption{Number of vertices ($N_v$), simplices ($N_{\Delta}$), and maximum edge length $h$ for the $k^{\text{th}}$ refinement of the icosahedral discretization of the sphere.}
\label{ico_table}
\end{table}

\subsubsection{Initial Conditions and Error Norms}

The forthcoming numerical tests, with the exception of the moving vortices test, are designed such that the initial condition returns to itself at the final integration time. Consequently, the backward characteristic map at time $T$ is given by the identity map on the sphere. We assess the accuracy of the method using an approximation of the following quantities 

\begin{subequations}
\label{error_measures}
\begin{align}
\text{tracer ($L^{\infty}$) error} &= \frac{\norm{ \mathcal{X}_{[T,0]}^* \phi_0 - \phi(\cdot,T)}_{L^{\infty}(\mathcal{S})}}{\norm{\phi(\cdot, T)}_{L^{\infty}(\mathcal{S})}}  \label{sup-norm}\, ,
\\
\text{tracer ($L^{1}$) error} &= \frac{\norm{ \mathcal{X}_{[T,0]}^* \phi_0 - \phi(\cdot,T)}_{L^{1}(\mathcal{S})}}{\norm{\phi(\cdot, T)}_{L^{1}(\mathcal{S})}}  \label{l1-norm}\, ,
\\
\text{map-(i) error} &= \norm{\mathcal{X}^{(i)}_{[T,0]} - \text{id}^{(i)}_{\mathcal{S}}}_{L^{\infty}(\mathcal{S})} \, ,
\\
\text{density error} &= \norm{1 - J_{\mu}(\mathcal{X}_{[T,0]})}_{L^{\infty}(\mathcal{S})}
\end{align}
\end{subequations}
where $\text{id}^{(i)}_{\mathcal{S}}$ is the $i$-th component of the identity map and $\phi(\cdot,T)$ is the expected solution at the final integration time. The sup-norm is approximated as the maximum value over $10^6$ points sampled from a random uniform distribution over the domain and the $L^1$ norm is approximated using a simple averaging over the vertices in the triangulation as 

\begin{equation}
\norm{f}_{L^{1}(\mathcal{S})} \approx \sum_{T \in \Delta} \sum_{v_i \in T} f(v_i)|T|/3 \, .
\end{equation}We note that it is common to assess the accuracy of a numerical method for linear advection using a discrete analogue of \eqref{sup-norm} in the $\ell^{\infty}$ norm. We have chosen a finer approximation of the continuous error measures \eqref{error_measures} here in an effort to more closely support the estimates given in section \ref{sec:error_estimates} and to emphasize the functional definition of the backward characteristic map. \par
We consider three different initial conditions; the first of which consists of two symmetrically located cosine-bells

\begin{equation}
g_i(\lambda, \theta) = 
\begin{cases}
 \frac{1}{2}\left[1 + \cos(\pi r_i/r)\right]  & \text{if  } r_i < r,  
\\
 0  \,\, & \text{otherwise, }
\end{cases}
\end{equation}
where $r = 1/2$ is taken to be the base radius of each cosine bell and $r_i = r_i(\lambda, \theta)$ is the great-circle distance from the centre of the bell located at $(\lambda_i, \theta_i)$, given by

\begin{equation}
r_i(\lambda, \theta) = \arccos(\cos \theta_i \cos \theta + \sin \theta_i \sin\theta \cos(\lambda - \lambda_i)) \, .
\end{equation}
The initial condition is then defined as

\begin{equation}
\label{cosine_bell}
\phi(R, \theta) = 0.1 + 0.9(g_1(R, \theta) + g_2(R,\theta)), 
\end{equation}and the centres of each bell are chosen to be $(\lambda_1, \theta_1) = (7\pi/6, \pi/2)$ and $(\lambda_2, \theta_2) = (5\pi/6, \pi/2)$. We note that for expressions involved spherical coordinates, we use the convention defined by the chart \eqref{spherical_coords}. The second initial condition, designed to assess the shape-preserving properties of the method, consists of two Zalesak disks \cite{zalesak1979fully} defined by

\begin{equation}
\label{zalesak_disk}
\phi(\lambda, \theta) = 
\begin{cases}
 1  & \text{if  } r_i \leq r \text{ and } |\lambda - \lambda_i| \geq r/6 \text{ for } i = 1,2 \, ,
\\
1  & \text{if  } r_i \leq r \text{ and } |\lambda - \lambda_1| < r/6 \text{ and } \theta - \theta_1 < -\frac{5}{12}r \, , \\
1  & \text{if  } r_i \leq r \text{ and } |\lambda - \lambda_2| < r/6 \text{ and } \theta - \theta_2 > \frac{5}{12}r \, ,\\
 0.1  \,\, & \text{ otherwise.}
\end{cases}
\end{equation}The third initial condition we consider is a sum of spherical harmonics of all degrees $(\ell, m)$ up to $ \ell = 32$ which are multiplied by a randomly generated real coefficients sampled from a uniform distribution from $-1$ to $1$. This initial condition is proposed by the authors since it is highly irregular and does not have compact support like the initial conditions \eqref{cosine_bell} and \eqref{zalesak_disk} do, providing a better assessment of the pointwise error introduced by the composition with $\mathcal{X}_{[t,0]}$. Based on the regularity properties of these initial conditions and the estimates provided in section \ref{sec:error_estimates}, the accuracy of the solutions for the cosine-bell (CB) and random spherical harmonics (rSph) are assessed using \eqref{sup-norm} whereas the solution for the slotted cylinder (SC) is assessed using \eqref{l1-norm}. We note that in every test performed, the approximation of $\mathcal{X}_{[T,0]}$ was accurate enough at the grid points such that the error \eqref{l1-norm} computed for the Zalesak disks was identically zero at every refinement except for $k = 0,1$, we have thus not included these results in the forthcoming convergence plots.  \par

In an effort to demonstrate the geometric flexibility of the method, the tests are performed on a sphere rotating about various axes. The rotating and inertial frames of reference are transformed between one another via a rotation matrix $R(t) \in \text{SO}(3)$ describing the rigid body rotation of the sphere. Let $\{\bsym{e}'_1, \bsym{e}'_2, \bsym{e}'_3\}$ be fixed basis vectors in the rotating frame such that $R(t)\bsym{e}'_i = \bsym{e}_i$ where $\bsym{e}_i$ is the standard Cartesian basis vector in the $i$-th coordinate direction in the inertial frame. If the velocity field has coordinate functions  $(u^1, u^2, u^3)$ in the rotating frame, then it can be expressed in the inertial frame as

\begin{equation}
\label{vel_transform}
\bsym{u}(\bsym{x},t) = \sum_{i=1}^3 u^i( R^T(t)\bsym{x},t) R(t) \bsym{e}'_i \,.
\end{equation}

\begin{figure}[!h]
\centering
\includegraphics[width = 15cm , height = 5cm]{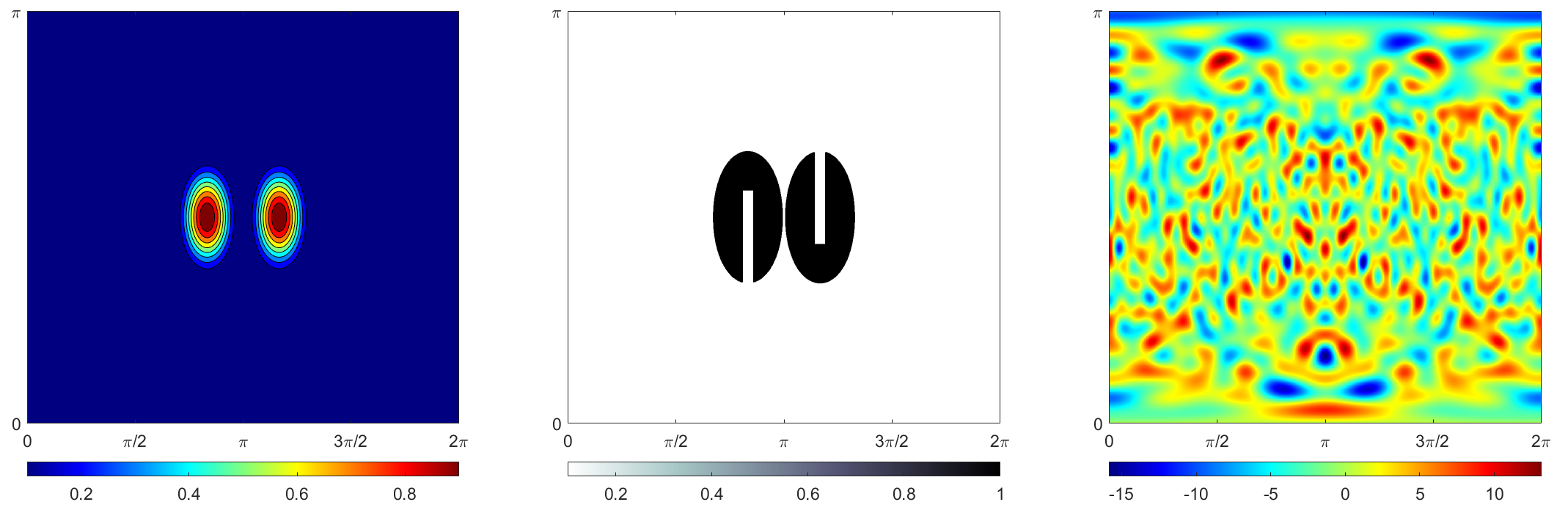}
\caption{Cosine bell \eqref{cosine_bell} (left), two Zalesak disks \eqref{zalesak_disk} (middle), and random spherical harmonics (right) initial conditions.}
\end{figure}

\subsubsection{Test case 1: Solid body rotation}
\label{sec:test_sbr}
\begin{figure}[h!]
\centering
\includegraphics[width = \linewidth , height = 6cm]{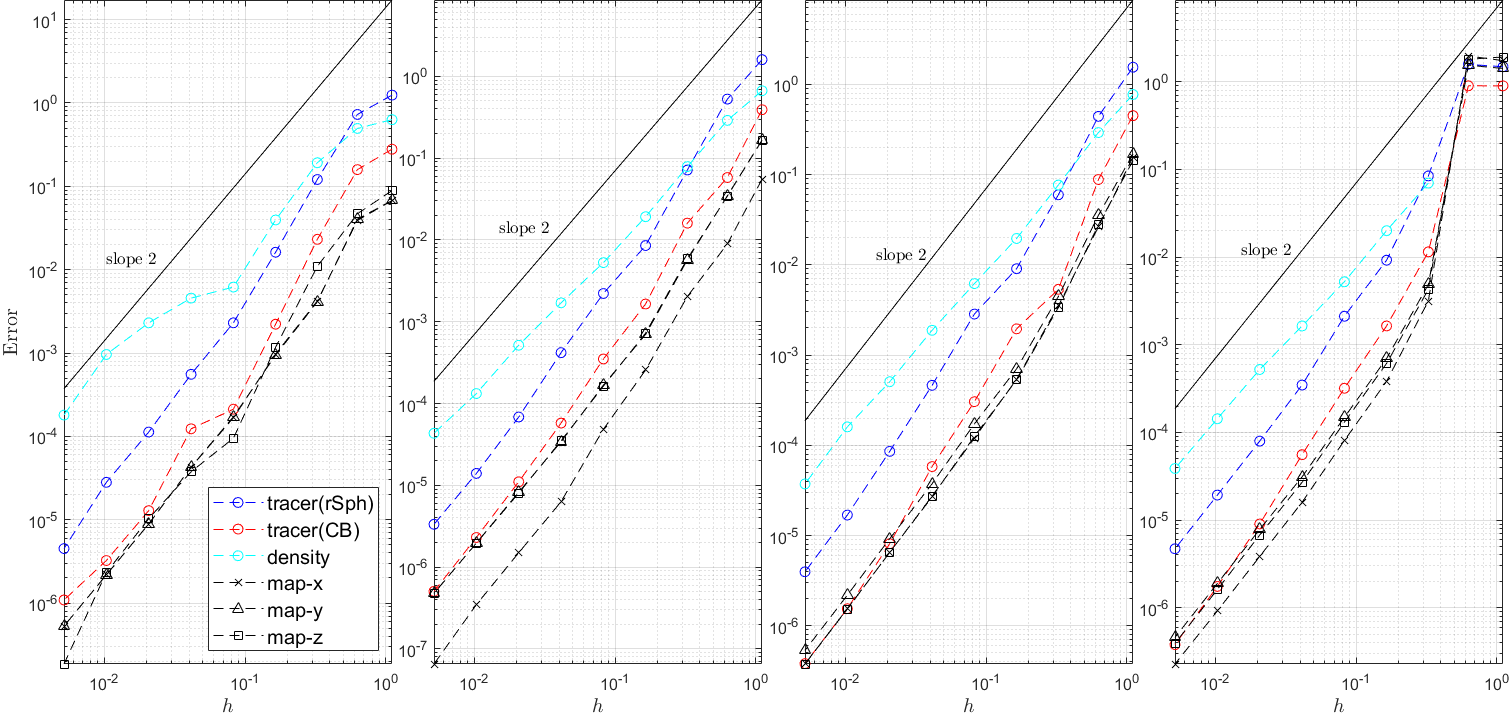}
\caption{Results for test case \ref{sec:test_sbr} for axes of rotation defined by $\alpha = 0, \pi/2, \pi/4, 1.05$ from left to right. The final integration time was taken to be $T=1$.}
\label{fig:sbr1} 
\end{figure}

The first test case considered is the solid body rotation for a rotation axis controlled by an angle of inclination $\alpha$ from the point $(0,0,1)$ along the direction of the positive $x$-axis. The velocity field is given by

\begin{equation}
\label{sbr_velocity}
\bsym{u}_{sbr} = \bsym{x} \times \nabla \psi \, , \quad \psi(x,y,z) = -\frac{2\pi}{T}(\cos\alpha z +\sin\alpha x)  \,,
\end{equation}
restricted to the sphere. The results of the test for values of $\alpha \in \{0, \pi/2, \pi/4, 1.05\}$ are shown in figure \ref{fig:sbr1}. We observe the order of convergence to be slightly greater than the expected second order convergence for each test case.

\subsubsection{Test case 2: Deformational flows}
\label{sec:deform_flow}
In this test case we consider five different deformational velocity fields. The first three velocity fields are the reversing deformational test combined with a solid body rotation about the axes defined by $\alpha = 0, 1.05, \pi/4$ \cite{nair2010class}. The transformation into the co-rotating frame is formed by rotation about the $y$-axis by $\alpha$ followed by a rotation about the $z'$- axis by angle $\beta(t) = -2\pi t/T$.  The velocity field is given by  

\begin{equation}
\label{bg_flow}
\begin{aligned}
\bsym{u}(\bsym{x},t) &= \bsym{u}_{sbr}(\bsym{x}) + \bsym{u}_{d}(\bsym{x},t) \,, 
\\
\bsym{u}_d(\bsym{x}', t) &= \bsym{x}'\times \nabla \psi_d(\bsym{x}',t) \, ,
\\
\psi_d(\bsym{x}',t) &= 2y'^2\cos(\pi t/T) \, \,,
\end{aligned}
\end{equation}restricted to the sphere. The gradient of the stream function is computed with respect to the primed coordinates and $\bsym{u}(\bsym{x},t)$ is evaluated in the inertial frame of reference using \eqref{vel_transform}. \par
The fourth and fifth velocity fields considered are the static and moving vortices \cite{nair2008moving}. The test consists of two circular vortices with antipodal centres positioned at the equator which we consider as rotating about the $z$-axis in the moving case.  Let $(\lambda',\theta')$ be spherical coordinates resulting from the transformation $\eta \circ \varphi_{\pi/2}(\bsym{x})$ where $\varphi_{\pi/2}$ is a rotation about the $x$-axis by angle $\pi/2$. The radial distance from the centre of the vortex is $\rho = \rho_0 \sin(\theta')$ and the angular velocity in dimensionless units is given by

\begin{equation}
\omega(\theta') = 
\begin{cases}
\frac{2\pi}{T} \frac{3\sqrt{3}}{2 \rho}\text{sech}^2(\rho) \tanh(\rho) \,\,\, &\text{if    } \rho \neq 0\, ,
\\
0  \,\,\, & \text{if    }  \rho = 0.
\end{cases}
\end{equation}
In the rotating frame of reference, we have that
$\dot{\lambda}' = \omega(\theta')\csc(\theta')$ , $\dot \theta' = 0 $ and the velocity field in the inertial frame of reference is given by 

\begin{equation}
\label{u_stat_vort}
\begin{aligned}
\bsym{u}(\bsym{x}) &= \bsym{u}_{s}(\bsym{x}) + \bsym{u}_v(\bsym{x}) \,, 
\\
\bsym{u}_{v}(\bsym{x}') &= -y' \omega(\theta')\bsym{e}'_1 + x'\omega(\theta')\bsym{e}'_2 \, ,
\end{aligned}
\end{equation}where $\bsym{u}_s$ is given by \eqref{sbr_velocity} for $\alpha = 0$ and the velocity field $\bsym{u}_v$ is evaluated in the inertial frame using \eqref{vel_transform} where $\varphi$ is a rotation about the $z$-axis by angle $\beta(t) = 2\pi t/T$. 
The solution to the advection equation for the velocity field \eqref{u_stat_vort} is given by

\begin{equation}
\label{vortex_solution}
\phi(\lambda', \theta', t) = 1 - \tanh\left[\frac{\rho}{5} \sin(\lambda' - \omega(\theta')t)\right],
\end{equation}
where the coefficients $1/5$ along with $\rho_0 = 3$ are chosen so that the deformation in the flow is smooth \cite{nair1999cascade, nair2002mass}. We test using both the rotating and static form of the vortex flow. \par

\begin{figure}[h!]
\centering
\includegraphics[width = 14cm, height = 4cm]{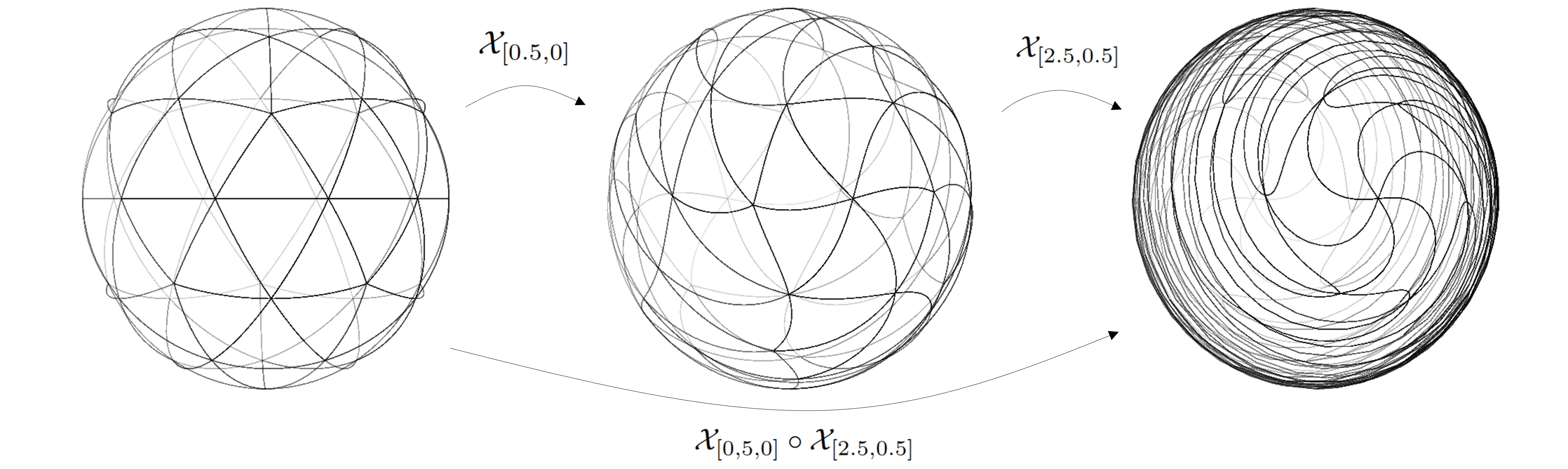}
\caption{The image of the first refinement of the icosahedral triangulation $\mathcal{X}_{[t,0]}(\mathcal{T})$ under the numerically approximated inverse flow map computed from \eqref{bg_flow} without rotation.}
\label{fig:deformed_mesh} 
\end{figure}

\begin{figure}[h!]
\centering
\includegraphics[width = \linewidth , height = 8cm]{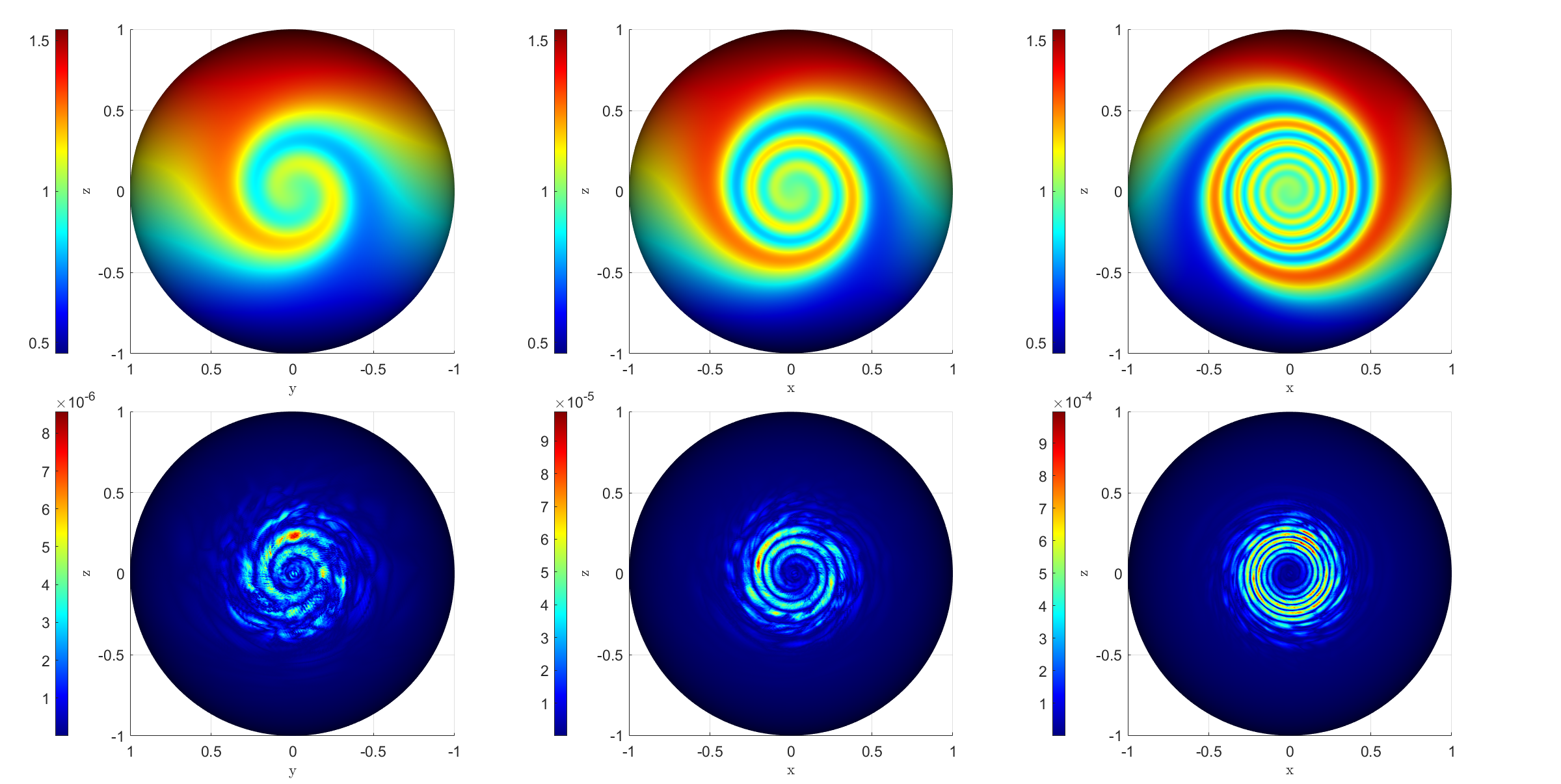}
\caption{Top row: Numerical solution for the moving vortices flow \eqref{u_stat_vort} for parameters $k = 8$ and $\Delta t \approx 0.0075$ at times $t = 0.5, 1, 2$ from left to right. Bottom row: The corresponding absolute value of the error from the solution \eqref{vortex_solution}.}
\label{fig:StatVort_Pics} 
\end{figure}

\begin{figure}[h!]
\centering
\includegraphics[width = \linewidth, height = 10cm]{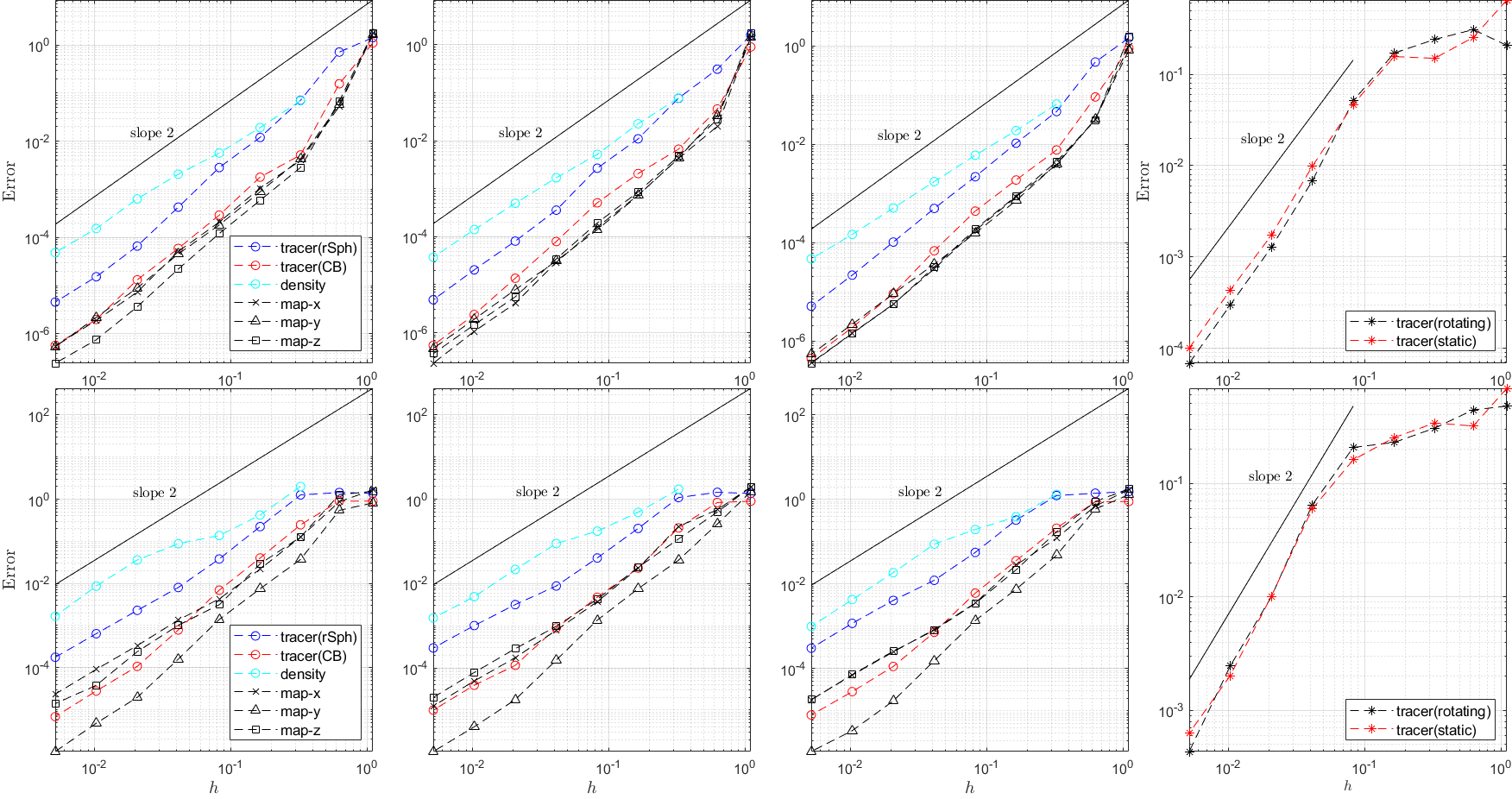}
\caption{Results for test case \ref{sec:deform_flow}, first three columns from left to right correspond to the velocity field \eqref{bg_flow} with background rotations defined by $\alpha = 0, \pi/4, 1.05$ respectively. The final integration times were taken to be $T=1$ (top row) and $T = 5$ (bottom row). The right most column corresponds to the velocity field \eqref{u_stat_vort} with and without rotation for final integration times $T=1$ (top) and $T = 2$ (bottom).}
\label{fig:deformational_flow_convergence_ico_T1} 
\end{figure}

The results of the numerical test are shown in figure \ref{fig:deformational_flow_convergence_ico_T1}. The tests are performed with final integration times of $T=1,5$ for the velocity fields defined by\eqref{bg_flow} and $T = 1,2$ for the vortices flow \eqref{u_stat_vort}. Increasing the final integration increases the error by approximately one order of magnitude, which is consistent with the error estimates and is attributed to the larger $T$ and $\Delta t$. The effect of the axis of rotation is insignificant for the velocity fields defined by \eqref{bg_flow} and we observe the theoretically expected second-order convergence. The effect of rotation is negligible for \eqref{u_stat_vort} and we observe the expected global second-order convergence, although this test required more refinements before reaching the asymptotic regime. We have included a visualization of the numerical solution to the moving vortices flow along with the error for the refinement $k = 8$, sampled at $1000^2$ points, in figure \ref{fig:StatVort_Pics}. The error is seen to become more localized about the centre of the vortex as time progresses.

\subsubsection{Test case 3: Compressible velocity field}
\label{sec:test_compressible}
In this test we consider a velocity field with non-zero divergence. Expressed in spherical coordinates, the velocity field is given by \cite{nair2010class} 

\begin{equation}
\label{div_flow_field}
\begin{aligned}
u(\lambda, \theta, t) &= -\sin^2(\lambda/2)\sin(2\theta)\sin^2(\theta)\cos(\pi t/T) \, ,
\\
v(\lambda, \theta, t) &= \frac{1}{2} \sin(\lambda)\sin^3(\theta) \cos(\pi t/T) \,  \, ,
\end{aligned}
\end{equation}which is transformed into Cartesian coordinates using \eqref{spherical_coords} and the fact that $(u,v) = (\dot{\lambda}(t)\csc(\theta(t)), \dot{\theta}(t))$. The test is designed to return to itself at the final integration time.\par

\begin{figure}[h!]
\centering
\includegraphics[width = 10cm, height = 6cm]{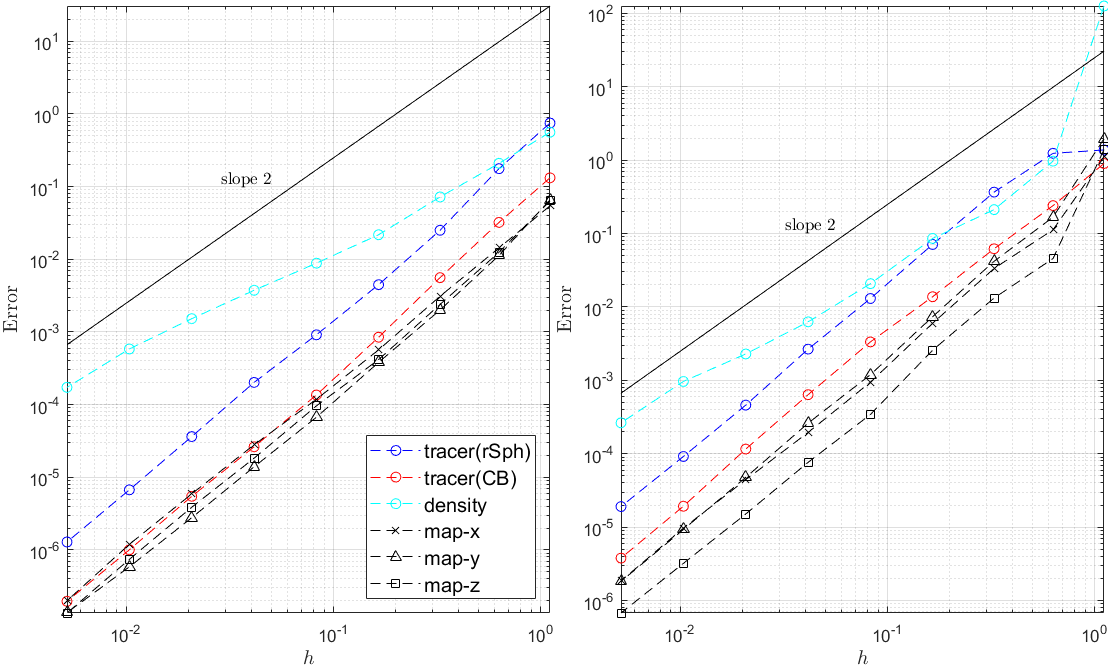}
\caption{Results for test case \ref{sec:test_compressible} with final integration times $T=1$ (left) and $T=5$ (right). }
\label{fig:test_divflow} 
\end{figure}

\begin{figure}[h!]
\centering
\includegraphics[width = \linewidth, height = 9cm]{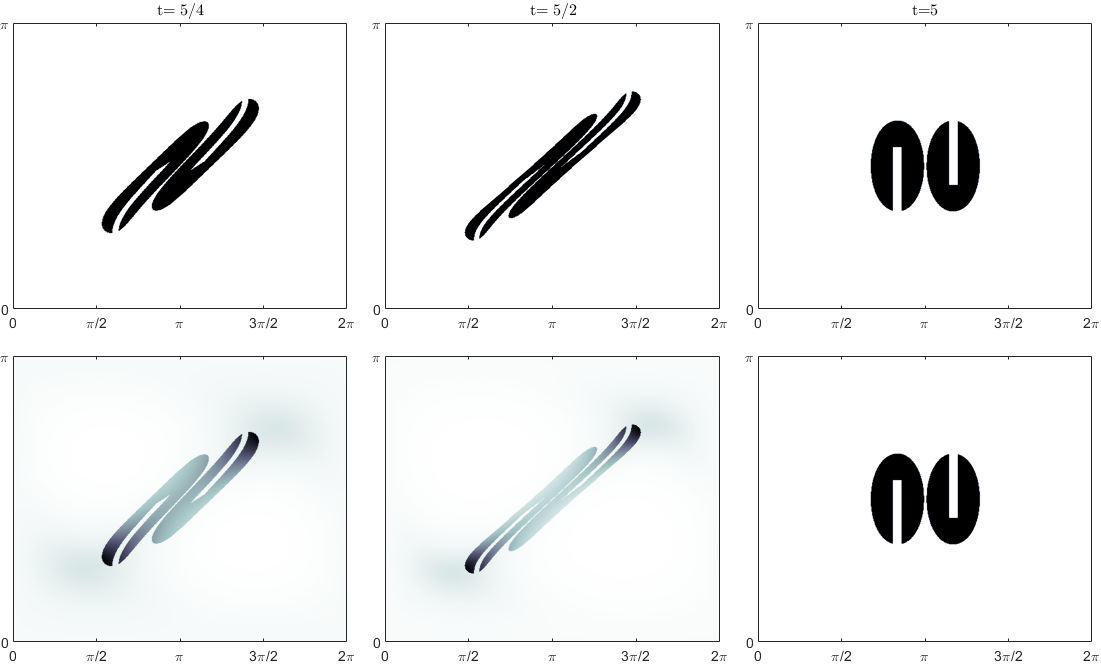}
\caption{Numerical solution for test case \ref{sec:test_compressible} with final integration time $T=5$ with initial condition \eqref{zalesak_disk}. Transported tracer $\phi_0 \circ \mathcal{X}_{[t,0]}$ (top row) and tracer density $\phi_0 \circ \mathcal{X}_{[t,0]}\cdot J_{\mu}(\mathcal{X}_{[t,0]})$ (bottom row) at times $t = 5/4, 5/2, 5$ from left to right.}
\label{fig:viz_divflow} 
\end{figure}

The compressibility of the velocity field manifests itself as a change in density of the fluid throughout the evolution. The numerical solution for the tracer and tracer density shown in figure \ref{fig:viz_divflow} illustrates this property. Since the action of the solution operator defined by \eqref{mass_solution_operator} and \eqref{solution_operator} are independent of the compressibility of the velocity field, the formulation needs no modification. The results of the convergence test are shown in figure \ref{fig:test_divflow}, affirming the theoretically predicted rate of convergence. 

\subsubsection{Test case 4: Correlated tracer transport}

It is important for transport schemes to preserve certain functional relations in the advected quantities. We evaluate the distortion of these relations due to the error introduced during computation using the diagnostic test case 3 provided in \cite{lauritzen2012evaluating}. As initial condition we consider two non-linearly correlated quantities, $q_1$ and $q_2$, defined by

\begin{equation}
\label{correlated_tracers}
\begin{aligned}
q_1(\lambda, \theta) &= 
\begin{cases}
 0.1 + 0.45\left[1 + \cos(\pi r_1/r)\right]  & \text{if  } r_1 < r,  
\\
 0.1 + 0.45\left[1 + \cos(\pi r_2/r)\right]  & \text{if  } r_2 < r,  
\\
 0.1  \,\, & \text{otherwise, }
\end{cases}
\\
q_2(\lambda, \theta) &= -0.8 q_1^2(\lambda,\theta) + 0.9 \, .
\end{aligned}
\end{equation}\begin{figure}[h!]
\centering
\includegraphics[width = 12cm, height = 5.5cm]{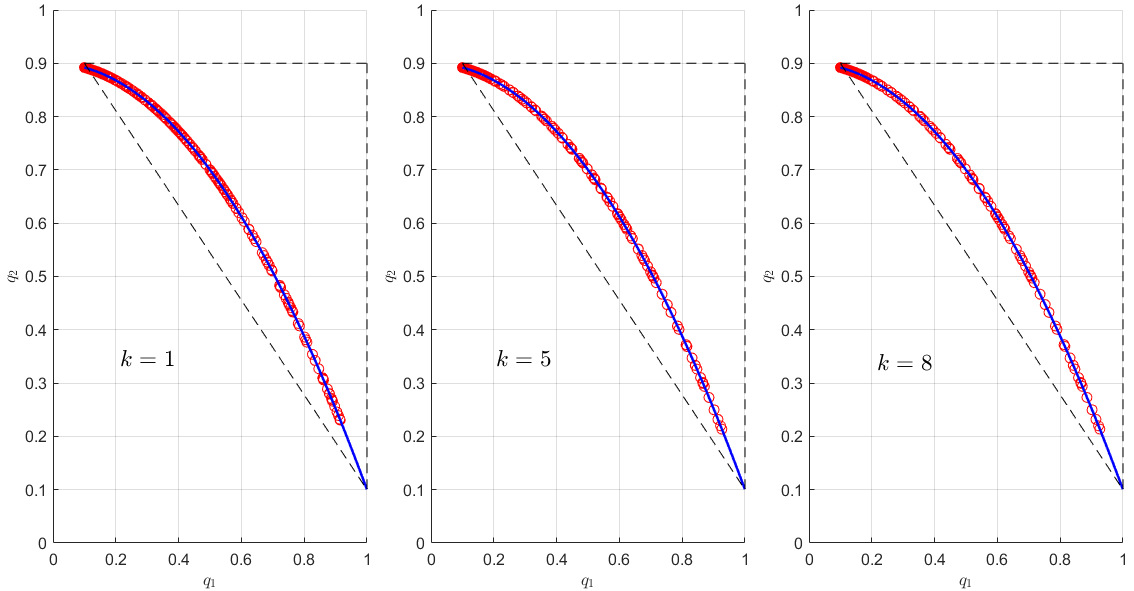}
\caption{Numerical mixing test for non-linearly correlated cosine-bell initial conditions \eqref{correlated_tracers} at refinements $k = 1,5,8$ (left to right), $q_2 \circ \mathcal{X}_{[T/2,0]}$ plotted against $q_1 \circ \mathcal{X}_{[T/2,0]}$, with grid spacings indicated. The transporting velocity field was \eqref{bg_flow}, final integration time $T =5$, initial conditions sampled at $t = 5/2$. The exact relation is indicated by the blue line and the computed results are plotted in red.}
\label{fig:numerical_mixing}
\end{figure}

As transporting velocity field, we use \eqref{bg_flow} with $\alpha = 1.05$ and final integration time $T=5$. The functional relation between these initial conditions is assessed using a scatter plot of $q_2$ plotted against $q_1$ at time $t=T/2$ when the deformation is the greatest. Each advected quantity is sampled at $200^2$ evenly spaced points in the computational domain for three different grid spacings. The results of the numerical test are presented in figure \ref{fig:numerical_mixing}. As a consequence of the problem formulation using the CM method and the analytic definition of the initial conditions, the functional relation between the tracers is preserved \emph{exactly}. Sampling each advected quantity at time $t = T/2$ via \eqref{solution_operator} requires a single evaluation of $\mathcal{X}_{[T/2,0]}$. As a result, the tracers are evaluated at the same location, thus exactly maintaining their correlations over time. \par 
We emphasize that the error incurred by the transported tracers is \emph{advective} in nature, resulting from an approximation of the location of evaluation. In this formulation, the effects of numerical mixing are not present since a single evaluation of the approximated backward characteristic map is used to advect each tracer. This property is not limited to two tracers and these results can be extended to many tracers with varying correlations while only incurring the added computational cost of the evaluation \eqref{solution_operator}.

\subsubsection{Test case 5: Mass conservation}
\label{conservation_test}

\begin{figure}[h!]
\centering
\includegraphics[width = \linewidth, height = 6cm]{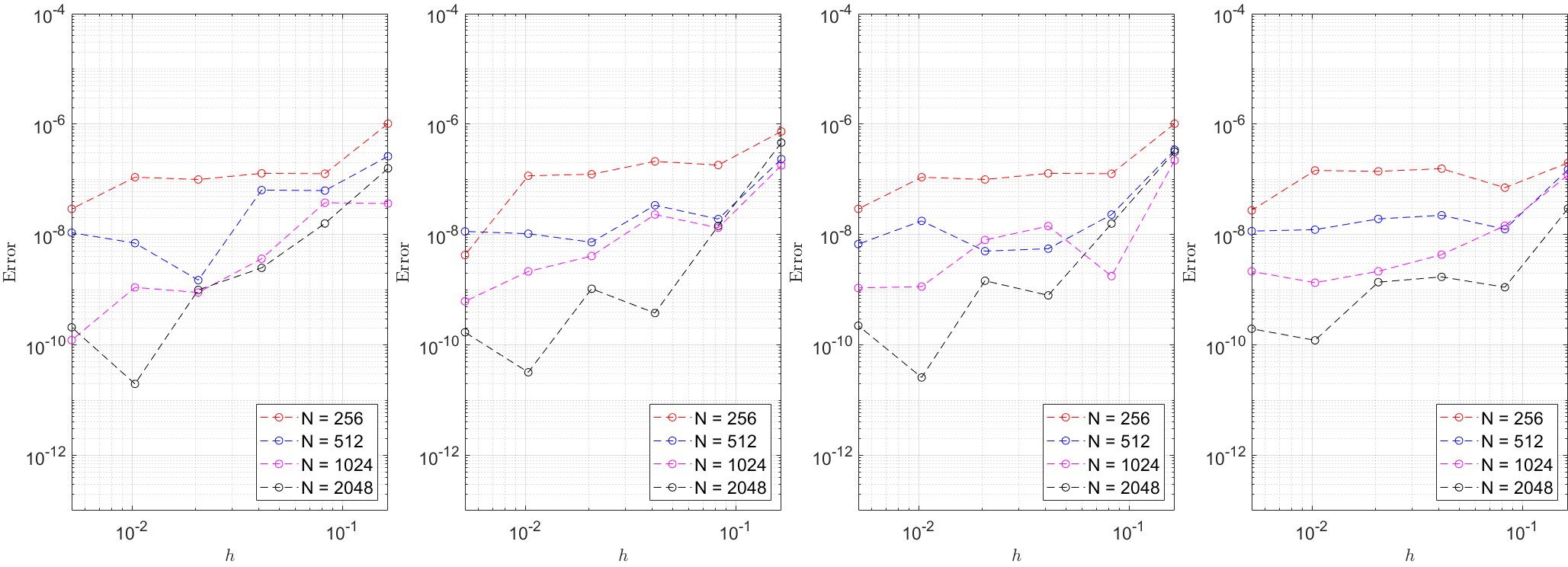}
\caption{Mass conservation errors at $t = T$ for test cases defined by the velocity field \eqref{bg_flow} with $\alpha = \pi/2, 1.05, \pi/4$ and the compressible velocity field \eqref{div_flow_field} from left to right. The parameter $N$ defines the number of grid points along each axis of the integration domain.}
\label{fig:conservation_errors} 
\end{figure}

As discussed in section \ref{sec:conservation_properties}, the method is conservative as a consequence of the change of variables formula. Any effort to demonstrate this property numerically will be a function of the quadrature scheme used to evaluate the integral. Let the normalized mass integral be  
\begin{equation}
\label{mass_functional}
I[\mathcal{X}_{[t,0]}] =  \frac{1}{4\pi}\int_{\mathcal{S}}\mathcal{X}_{[T,0]}^*\mu  = \frac{1}{4\pi} \sum_{k = 1}^N \int_{V_k} \mathcal{X}_{[t,0]}^*\mu ,
\end{equation}where $V_k \subset \mathcal{S}$ and partition the sphere. We seek to demonstrate an independence of the mass conservation from the parameters $N_t$ and $N_{\Delta}$. We used a $9$ point Gauss-Legendre quadrature scheme to approximate the integrals over elements $V_k = \eta^{-1}(A_{k})$ where the $A_k$ form a rectangular partition of $(\lambda, \theta) \in [0,2\pi) \times [0,\pi]$ with $N$ grid points along each axis. We have considered the tests defined by the velocity fields \eqref{bg_flow} about various axes of rotation and \eqref{div_flow_field} up to the final integration time $T = 5$. The results are given in figure \ref{fig:conservation_errors} and demonstrate a high-degree of accuracy for the mass error at $t = T$ for all refinements. Although we were unable to absolve the error completely from a dependence on the mesh size $h$ for each test, it is observed that the mass error for meshes defined by more than three refinements overall decrease as we refine the numerical quadrature grid with respect to $N$. These results are consistent with the mass conserving property of the method and the strongest agreement was for the compressible test case, where the error is more apparently independent of the discretization of $\mathcal{X}_{[t,0]}$.

\subsection{Submap Decomposition}

\begin{figure}
    \centering
    \includegraphics[width = 12cm, height = 5.5cm]{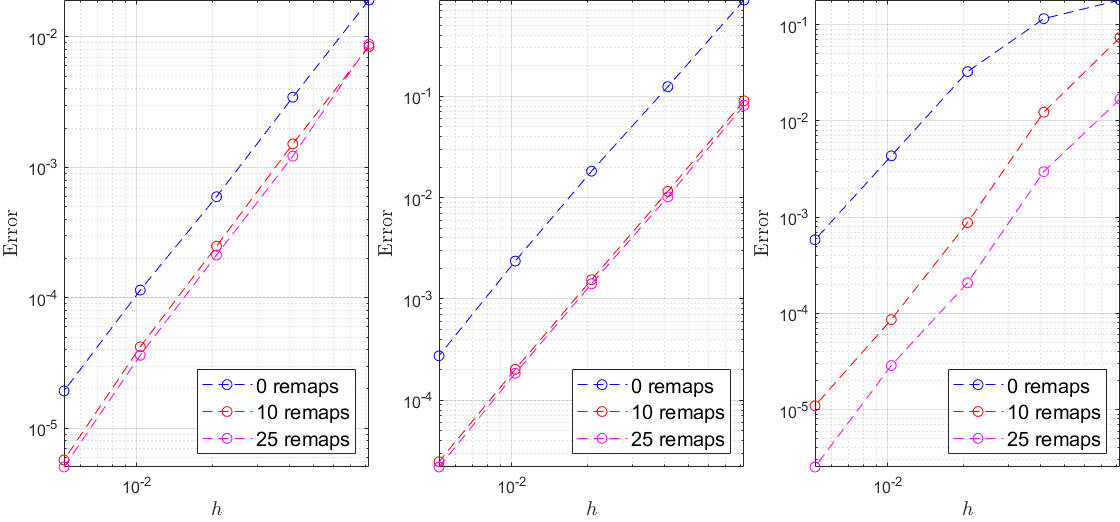}
    \caption{The error \eqref{sup-norm} for the random spherical harmonic tracer at different refinement levels and remapping strategies. We have considered the test cases defined by the velocity fields \eqref{u_stat_vort} for $T=2$, \eqref{bg_flow} for $\alpha = 1.05$ and $T=5$, along with \eqref{div_flow_field} for $T=5$.}
    \label{fig:remapping_convergence}
\end{figure}

In this section we present results for a fixed remapping technique and compare the accuracy gained by increasing the amount of remapping steps. We compare the error \eqref{sup-norm} for the test cases defined by the velocity fields \eqref{u_stat_vort} for $T=2$, \eqref{bg_flow} for $\alpha = 1.05$ and $T=5$, along with \eqref{div_flow_field} for $T=5$. The simulations were run using $N_t = 250$ time steps for $10$ and $25$ remaps. The results of the test are given in figure \ref{fig:remapping_convergence} where we have included the case without remapping for comparison. The accuracy is shown to improve as the number of remapping steps increases with the greatest improvements observed for the rotating vortices test case. It is observed that the accuracy does not increase significantly between the tests for $10$ remaps compared to $25$ remaps, suggesting a saturation of the error as the number of remaps increases. An analysis of the optimal number of remaps for a given velocity and integration time is warranted and the subject of our current research.     

\subsection{Fractal Set Advection} \label{last_section}
\begin{figure}[h!]
\centering
\begin{tabular}{cc}
\includegraphics[width = 0.4\linewidth, height = 5.7cm, trim = {2cm 2cm 0 1cm}, clip]{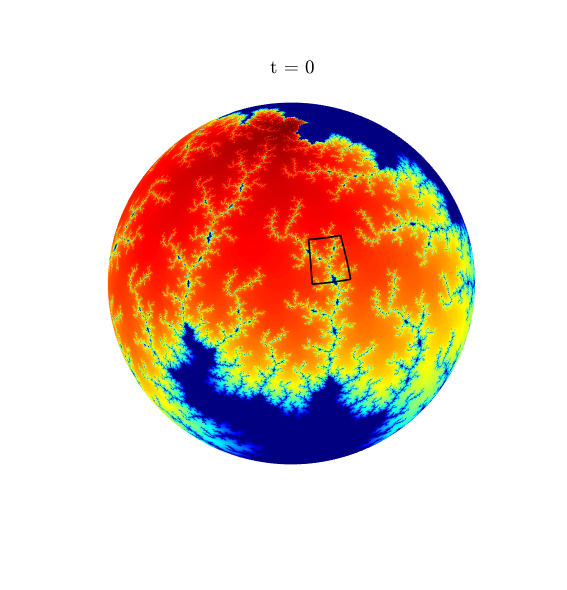} & \includegraphics[width = 0.4\linewidth, height = 5.7cm, trim = {2cm 2cm 0 1cm}, clip]{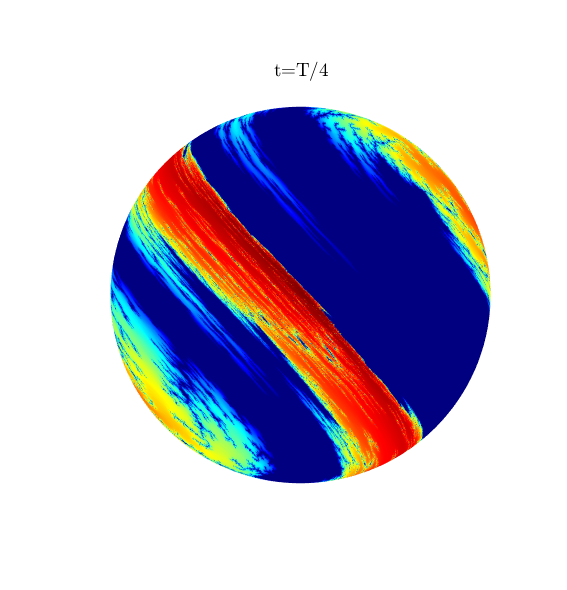} \\
\includegraphics[width = 0.4\linewidth, height = 5.7cm, trim = {2cm 2cm 0 1cm}, clip]{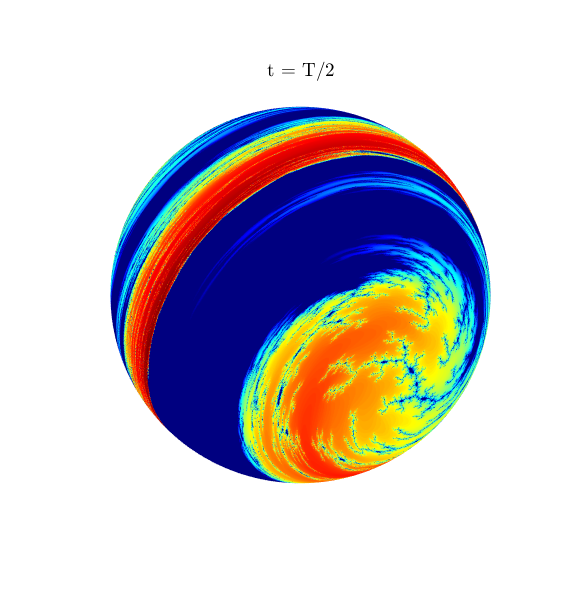} & \includegraphics[width = 0.4\linewidth, height = 5.7cm, trim = {2cm 2cm 0 1cm}, clip]{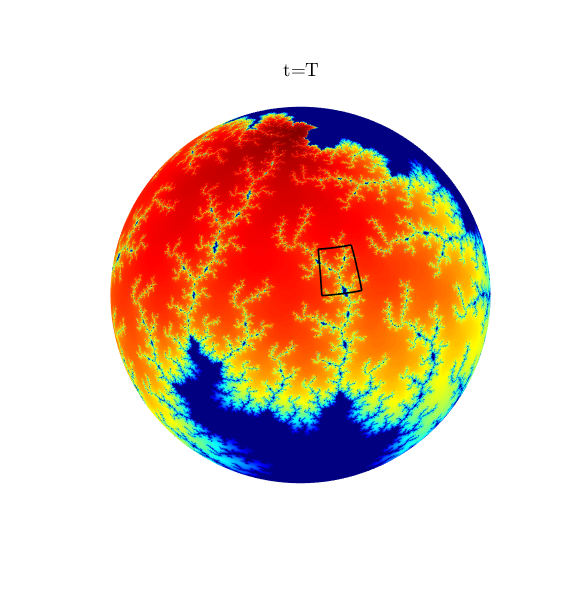} \\
\end{tabular}

 \includegraphics[width = 14cm, height = 7.2cm]{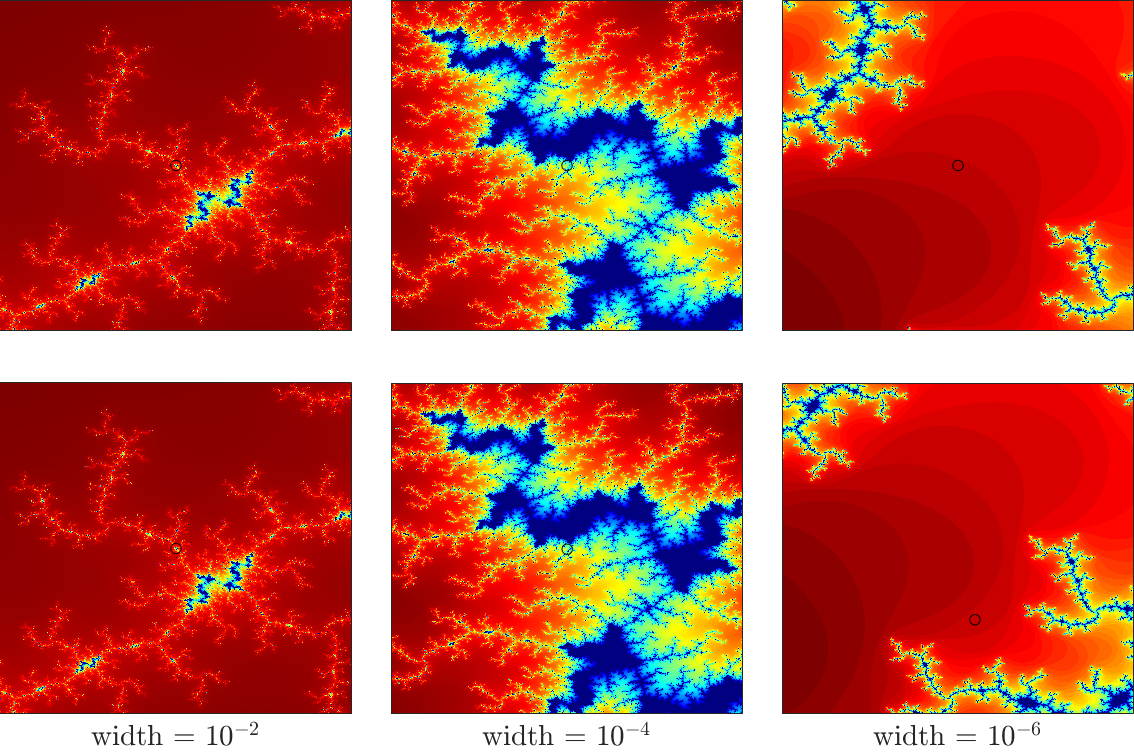}

\footnotesize
 \caption{Top: Advection of the Mandelbrot set on the sphere under the velocity field \eqref{bg_flow} with parameters $\alpha = \pi/4$ and $T = 5$. A window of angular width $2^{-2}$ centred at the focal point of the zoom is included. Bottom: Zoom on window depicted at times $t = 0$ (top row) and $t = 5$ (bottom row) up to a width of $10^{-6}$ with focal point given by the black circle.}
 \label{fig:Mandelbrot}
 \end{figure}

We present a numerical experiment designed to illustrate the capacity of the method to advect functions with poor regularity and represent fine features by considering the advection of the Mandelbrot set on the sphere. The set is generated via stereographic projection with a cap of arclength $10^{-9}$ at the north pole. The origin of the stereographic plane (i.e. the south pole) is placed at the point $-0.235125 + 0.827215i$ in the complex plane defining the Mandelbrot set. The axes are then scaled by a factor of $4\times 10^{-5}$. The velocity field is taken to be \eqref{bg_flow} with $\alpha = \pi/4$ and $T=5$ and we advect the set using the backward characteristic map $k=8$ refinement in the convergence test presented in figure \ref{fig:deformational_flow_convergence_ico_T1}. We demonstrate a gradual zoom to a frame size of $10^{-6}$ in figure \ref{fig:Mandelbrot} for the solution at times $t = 0$ and $t = 5$. Each image is produced by sampling the map at $1600^2$ points within the frame. This final frame size is chosen based on the results of figure \ref{fig:deformational_flow_convergence_ico_T1}, where we expect to observe a discrepancy between the initial condition and the solution at the final time for a window width on the same order of magnitude as the $L^{\infty}$ error of the backward characteristic map, as observed in the final column of the zoom \ref{fig:Mandelbrot}. 

\setlength{\tabcolsep}{0pt}

\section{Conclusion}

The CM method for the numerical solution of the tracer transport equations on the sphere has been presented. The method employs an extrinsic approach to the spatial discretization of the backward characteristic map using a projection-based interpolation framework for manifold-valued data. The resulting globally defined approximation acts via pullback as a solution operator to the tracer transport equations on the sphere. Moreover, the method is formulated to be independent of the compressibility of the transporting velocity field. The error estimates provided in section \ref{sec:error_estimates} are supported in section \ref{convergence_tests} by convergence tests comprised of solid body rotations about various axes, deformational flows, and compressible flows. The method is shown to be globally second-order accurate for tracer transport independent of the velocity and second-order for the density transport in the incompressible setting and first-order accurate for compressible flows. In addition to these results, we demonstrated some properties of the advective nature of the error introduced by the method with the ability to exactly preserve non-linear correlations between tracers, the conservation of mass, and the transport of a fractal set in a complex flow environment.

Improving upon the method in its present form offers many interesting avenues of investigation. Given its unique features, extending the method to incorporate reaction and diffusion terms for the advection of multiple chemical tracer species and the comparison with existing techniques is warranted. An operator splitting technique where the method presented here is used for the transport operator may be a viable strategy for this purpose. Beyond transport on the sphere, the extrinsic formulation presented here permits the application of the method to more general embedded submanifolds. Finally, the method is not limited to the transport of quantities which do not interact with the advecting velocity field. The CM method has been successfully applied to the solution of the incompressible Euler equations on a two dimensional periodic domain using the advection of the vorticity \cite{yin2021characteristic}. The extension of these techniques to a spherical geometry and to the solution of the shallow water equations is the subject of our current research. We believe that the CM method offers a unique and useful tool for problems related to geophysical fluid dynamics.

\section*{Declaration of competing interest}

The authors declare that they have no known competing financial interests or personal relationships that could have appeared to influence the work reported in this paper.

\section*{Acknowledgments}

The work of S.T. was partially supported by the NSERC CGS-D program.  The work of J-C.N. was partially supported by the NSERC Discovery Grant program and the Agence Nationale de la Recherche (ANR), grant ANR-20-CE46-0010-01. The authors would like to thank the anonymous referees whose constructive comments have improved the article and Xi-Yuan Yin of McGill University for many helpful discussions and insights.

\section*{Appendix: Construction of the Powell-Sabin Spherical Spline}
\label{appendix}
In an effort to better facilitate the reproducibility of the method presented, we provide details on the explicit construction of the Powell-Sabin spherical spline macro-element method, following the descriptions given in \cite{lai2007spline, schumaker2015spline, alfeld1996fitting}. In a spherical triangle $T \in \mathcal{T}$ with vertices $v_1, v_2, v_3 \in \mathbb{R}^3$, a local coordinate system is defined using the spherical barycentric coordinates of a point $\bsym{v} \in T$. These are given by $\bsym{b}(\bsym{v}) = (b_1, b_2, b_3) \in \mathbb{R}^3$ such that $\bsym{v} = b_1v_1 + b_2v_2 + b_3v_3$. Nearly all properties of trihedral barycentric coordinates carry over to the spherical case, except that they are not required to add up to one \cite{alfeld1996bernstein}. A local spherical Bernstein-B{\'e}zier (SBB) polynomial in $T$ is defined by
\begin{equation}
\label{quad_HBB}
p(v) = \sum_{i+j+k = d} c_{ijk}B^{d}_{ijk}(v) \, , \quad B^{d}_{ijk}(v) = \frac{d!}{i!j!k!}b^{i}_1b^{j}_2b_3^k \, , \quad i + j+ k = d\, .
\end{equation}
for coefficients $c_{ijk} \in \mathbb{R}$. Efficient evaluation of a spherical polynomial in Berstein-B{\'e}zier form can be performed using deCasteljau's algorithm \cite{lai2007spline}. The directional derivative of a spherical polynomial $p$ in the direction $\bsym{g} \in \mathbb{R}^3$ at the point $\bsym{v} \in T$ is given by 
\begin{equation}
\label{SBB_derivative}
D_{\bsym{g}} p(\bsym{v}) = \bsym{b}(\bsym{g})\cdot\nabla_{b}p(\bsym{v}) \, ,
\end{equation}where $\bsym{b}(\bsym{g})$ are the barycentric coordinates of $\bsym{g}$ relative to $T$ and $\nabla_b$ is the derivative of the spherical polynomial in each barycentric coordinate \cite{alfeld1996fitting}.  We note that the directional derivative of a SBB polynomial is defined with respect to a homogeneous extension from the sphere, and is independent of this extension in the directions $\bsym{g}$ which are tangent to the sphere \cite{alfeld1996fitting}. \par

\begin{figure}[!h]
\centering
\includegraphics[width = 12cm , height = 5cm]{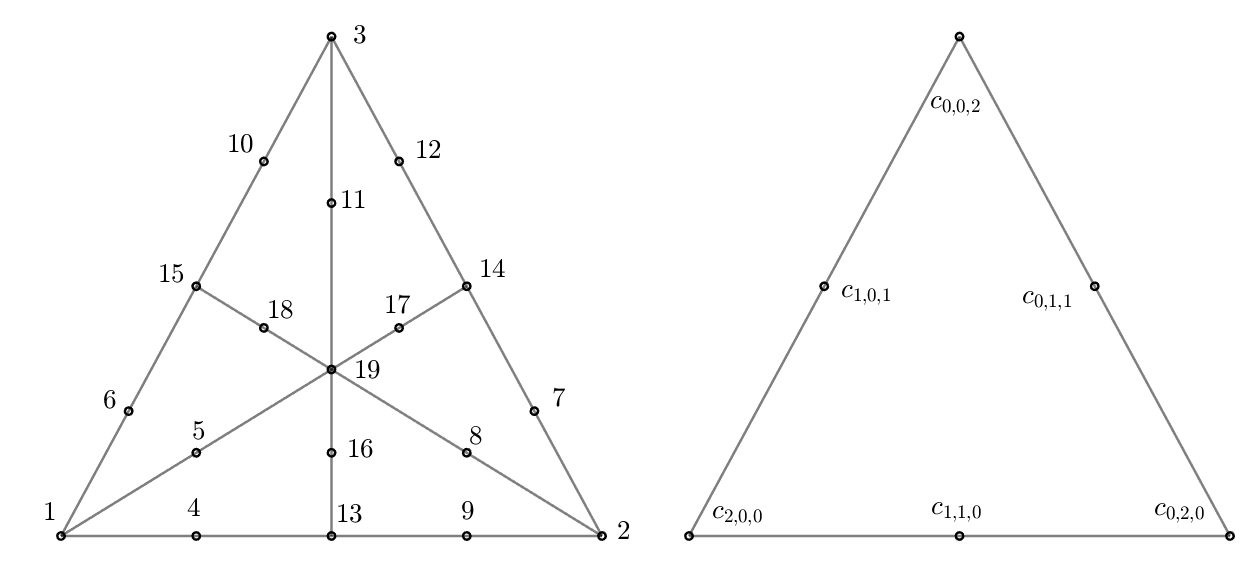}
\caption{Left: Powell-Sabin split and labelled coefficients. Right: Coefficients associated to the quadratic polynomial on the sub-triangles.}
\label{fig:ps_split} 
\end{figure}

The Powell-Sabin spherical spline defined by the Hermite interpolation operator \eqref{Hermite_operator} can be written explicitly by combining the $9$ pieces of data with the $C^1$ continuity conditions across the edges of the subdivision to map into the space $S_2^1(\mathcal{T}_{PS})$ as a set of $19$ coefficients for each macro-triangle $T$. The evaluation of a spline defined on the PS-split at a point $\bsym{v} \in \mathcal{S}$ is performed by first determining the macro-triangle $T$ through a containing simplex querying strategy. Then, using the barycentric coordinates of $\bsym{v}$ relative to $T$ one can determine in which subdivided triangle the point $v$ is contained. The barycentric coordinates of $\bsym{v}$ are then computed relative to the sub-triangle and a quadratic SBB polynomial of the form \eqref{quad_HBB} where the $c_{ijk}$ are given by the $6$ out of the $19$ coefficients which are associated with the sub-triangle (left side of \ref{fig:ps_split}), arranged appropriately as in the right panel \ref{fig:ps_split}. Let $v_4$ be the spherical barycenter of the triangle and denote $v_{ij}$ the midpoint of the edge connecting the vertex $v_i$ to $v_j$. Using the definition of the barycentric coordinates, it is readily observed that $c_1 = f(v_1), c_2 = f(v_2), c_3 = f(v_3)$. The coefficients $c_4, c_5, \dots, c_{12}$ can then be determined through the directional derivative formula \eqref{SBB_derivative}. In particular, suppose that $(a_1, a_2,0)$ and $(\tilde{a}_1,0, \tilde{a}_3)$ are the barycentric coordinates relative to the sub-triangle $T_1 = \langle v_1, v_{12},v_4\rangle$ of the unit vectors $\bsym{e}_{12}$ and $\bsym{e}_{14}$ tangent to the sphere at $v_1$ pointing toward $v_{12}$ and $v_4$ respectively. Using \eqref{SBB_derivative} and $\bsym{b}(v_1) = (1,0,0)$ we get $c_4, c_5$ as 

\begin{equation}
c_4 = \left(D_{\bsym{e}_{12}}s(v_1)/2 - a_1 c_1\right)/a_2 \,, \quad c_5 = \left(D_{\bsym{e}_{14}}s(v_1)/2 - \tilde{a}_1 c_1\right)/\tilde{a}_3 \,.
\end{equation}
The coefficients $c_6, \dots , c12$ can be determined in the other subtriangles analogously. We note that when computing these quantities it helps to arrange the six sub-triangles as $T_1 = \langle v_1, v_{12}, v_4 \rangle,\, T_2 = \langle v_2, v_{4}, v_{12} \rangle,\, T_3 = \langle v_4, v_{2}, v_{23} \rangle,\, T_4 = \langle v_{23}, v_{4}, v_{3} \rangle,\, T_5 = \langle v_{13}, v_4, v_{3} \rangle,\, T_6 = \langle v_4, v_{13}, v_1 \rangle$ where the sub-triangles are ordered counterclockwise from the first vertex. Computing the barycentric coordinates in each sub-triangle with these orientations neatly enforces continuity across the edges of the subdivision. The remaining coefficients on the edges are given by

\begin{equation}
\begin{aligned}
c_{13} &= r_1c_4 + s_1c_9\,,\quad c_{14} = r_2c_7 + s_2 c_{12}\,,\quad c15 = r_3c_{10}
+ s_3c_{6}\,, 
\\
c_{16} &= r_1c_5 + s_1c_8 \,,\quad c_{17} = r_2c_8 + s_2c_{11}\,,\quad c{18} = r_3c_{11} + s_3c_5 \,, 
\\
c_{19} &= a_1 c_5 + a_2 c_8 + a_3 c_{11} \,,
\end{aligned}
\end{equation}
where $(r_i,s_i)$ are given by the spherical barycentric coordinates of the macro-triangle edge midpoints such that $e_i = r_iv_i + s_iv_{i+1}$ for $i = 1,2,3$ (where $v_{3+1} = v_1$) and the $(a_1,a_2,a_3)$ are the spherical barycentric coordinates of $v_4$ with respect to the macro triangle.

\FloatBarrier
\footnotesize
\bibliographystyle{elsarticle-num} 
\bibliography{advection_S2}

\end{document}